\documentclass{amsart}
\usepackage{dcolumn}
\usepackage{graphicx}
\usepackage{amsmath}

\newtheorem{exa}{Example}
\newtheorem{defn}{Definition}
\newtheorem{remark}{Remark}

\usepackage{color}
\definecolor{brown}{rgb}{0.8,0.6,0.3}
\definecolor{dgreen}{rgb}{0.2,0.4,0.3}

\usepackage[pdfauthor={Richard J. Mathar}, colorlinks=true,urlcolor=blue,linkcolor=red,citecolor=brown,
pdfkeywords={{simple graph},{cycle},{tree},{Otter}},pdftitle={Connected Graphs with Isolated Cycles},bookmarksopen=true,pdfpagelayout=OneColumn]{hyperref}

\begin{document}

\title{Counting Connected Graphs without Overlapping Cycles}

\author{Richard J. Mathar}
\urladdr{http://www.mpia.de/~mathar}
\address{
Hoeschstr.\ 7, 52372 Kreuzau, Germany
}

\date{\today}

\begin{abstract}
The simple connected graphs may be classified by their cycle composition (number
and lengths of cycles). This work derives the counting series of the simple
connected graphs that have cycles of unrestricted number and length, but no
overlapping cycles. Cycle pairs of these graphs of interest
must not have common nodes or edges.

The recipe of counting
these graphs is based on the counting series of the associated planted graphs,
multisets of planted graphs, a recursive synthesis of enriched trees, and a
generalized Otter's formula that maps the underlying rooted block graphs to
the underlying block graphs.
\end{abstract}

\maketitle

\section{Graph Cycle Structures}
The simple connected graphs \cite[A001349]{EIS} may be classified
by their cycle structure. Classifications may impose restrictions on
the number of cycles, the lengths of cycles, the degrees of nodes,
or other criteria, see for example:
\begin{exa}
The subset of graphs without cycles, called trees \cite[A000055]{EIS};
\end{exa}
\begin{exa}
the unicyclic graphs with cycles of length $\ge 3$ \cite[A001429]{EIS};
\end{exa}
\begin{exa}
The connected graphs with at most one cycle, no cycle of length 2 \cite[A005703]{EIS};
\end{exa}
\begin{exa}
the unicyclic graphs with cycles of length $\ge 1$ and node
degrees not larger than 4 \cite[A002094]{EIS}.
\end{exa}

This work enumerates a set of graphs posing a restriction
on the entanglement of the cycles:
\begin{defn}
A C-tree is a simple (loopless, non-oriented, connected) graph where
cycle pairs do not overlap, which means any pair of two cycles in the
graph must not contain a common node.
\end{defn}
\begin{defn}
The ordinary generating function for C-trees on $n$ nodes is
\begin{equation}
C(x) = \sum_{n\ge 0} c(n)x^n = 1+x+2x^2+\cdots
\end{equation}
\end{defn}

\begin{remark}
C-trees require that any two cycles are not just sharing
no edge but also no node. In that respect C-trees are
a narrower concept than Husimi graphs \cite{HararyPNAS39,HararyAnM58}.
\end{remark}

There is no specification with regard to the length of the cycles
in a C-tree. The multiset of cycles may include cycles of length 2
(which will appear as double edges in the illustrations that follow).
For the combinatorial enumeration it will often by useful to regard
nodes that are not part of any cycle as cycles of length 1\@.
So the trees are C-trees where all cycles have length 1.

Because each node in a C-tree is a member of exactly one cycle
(allowing cycles of length 1), there is a surjection of the set
of C-trees to trees by contracting the nodes common to a cycle
to a single node:
\begin{defn}
The skeleton tree of a C-tree is the (ordinary) tree obtained
by collapsing the nodes of each cycle into a single node. The
edges of degree 2 that connect pairs of nodes within a cycle are deleted,
then the set of nodes in the (former) cycle is replaced
by a single node, pulling the edges that led to branches
outside the cycle into the single node. 
\end{defn}
\begin{remark}
The bridges of the C-tree
are the edges in the skeleton tree.
\end{remark}
\begin{remark}
The C-trees contain two types of blocks: the cycle graphs and the tree of two
points ($K_2$ graph). The skeleton tree is the block graph of the C-tree.
\end{remark}

\begin{figure}
\includegraphics[scale=0.4]{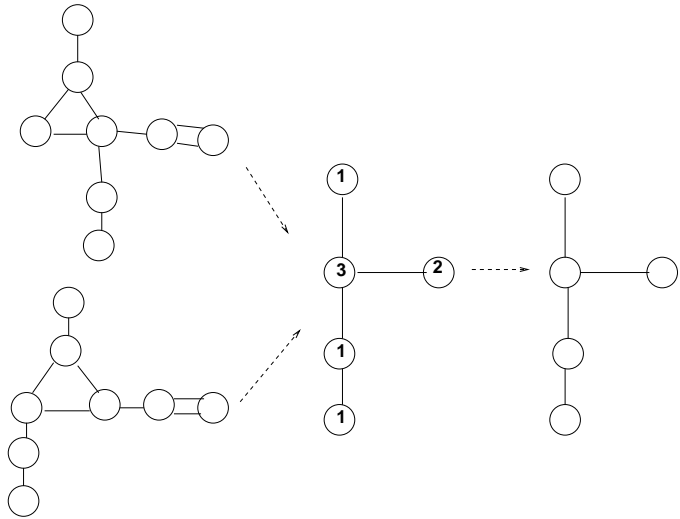}
\caption{Example of a surjection of two C-trees with 8 nodes
to a weighted tree of 5 nodes and further surjection to their skeleton
tree of 5 nodes.}
\label{fig.contr}
\end{figure}

The contraction of a C-tree to its skeleton may pass
through an intermediate step of a weighted tree, where
the weights are the lengths of the cycles of the C-tree;
see Figure \ref{fig.contr}.

\begin{remark}
Mapping C-trees to weighted trees loses the information
on to which nodes of the cycles the external edges
were attached.
So the number of weighted trees with weight $n$
\cite[A036250]{EIS} is a lower bound on $c(n)$.
\end{remark}

\begin{remark}
The C-trees are planar graphs.
\end{remark}

\section{Hand-counted C-trees Up to 6 Nodes}
\subsection{Injection Methology}

To count C-trees by hand, we interpret the number $n$ of nodes
as the weight of weighted trees \cite{HararyAM101,MatharVixra1805}, working 
backwards through the surjection if Figure \ref{fig.contr}.
We puff up each node of weight $c$ to a cycle of length $c$,
generating one or more C-trees.
\begin{defn}
An endnode of a graph is a node of degree 0 or 1, which means,
it has at most one edge attached to it.
\end{defn}
If the weight is $c>1$ and not on an endnode, more than one C-tree
may emerge.
\begin{defn}
$T(v)$ is the number of trees with $v$ nodes \cite[A000055]{EIS},
\begin{equation}
T(v)=1,1,1,2,3,6,11,23,47,106,235,\ldots \quad (v\ge 1).
\end{equation}
\end{defn}
\begin{defn}
$r_t(n)$ is the number of C-trees with $n$ nodes that
have skeleton trees with $t$ nodes:
\begin{equation}
c(n) = \sum_{t=1}^n r_t(n). 
\end{equation}
\end{defn}

In particular
\begin{equation}
r_1(n)=1
\end{equation}
because the skeleton tree with one node can be puffed up
in exactly one way to the cycle graph of $n$ nodes.
\begin{equation}
r_n(n)=T(n)
\end{equation}
is the number of trees with $n$ nodes, meaning that
the skeleton tree \emph{is} the C-tree.
If the skeleton tree has two nodes, the two weights of the
associated weighted trees may be any partition of $n$ into 2 parts:
\begin{equation}
r_2(n) = \lfloor n/2\rfloor.
\end{equation} 

\subsection{Zero or one Node}
\begin{figure}
\includegraphics[scale=0.3]{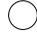}
\caption{The C-tree with 1 node.}
\label{fig.R1}
\end{figure}

There is one empty graph and one C-graph with a single node (Figure \ref{fig.R1}).
\begin{equation}
c(0)=c(1)=1.
\end{equation}
\subsection{Two nodes}
\begin{figure}
\includegraphics[scale=0.3]{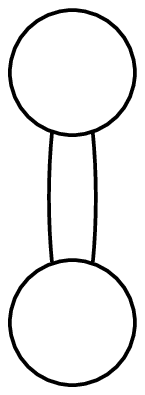}
\includegraphics[scale=0.3]{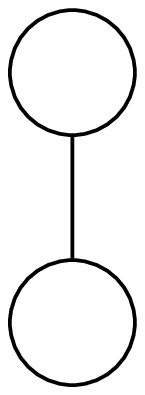}
\caption{The 2 C-trees with 2 nodes.}
\label{fig.R2}
\end{figure}

(i) The skeleton tree on one node generates
one C-tree with a cycle of length 2,
$r_1(2)=1$.
(ii) The skeleton tree on two nodes contributes that tree, $r_2(2)=1$
(Figure \ref{fig.R2}):
\begin{equation}
c(2)=1+1=2.
\end{equation}
\subsection{Three nodes}
\begin{figure}
\includegraphics[scale=0.3]{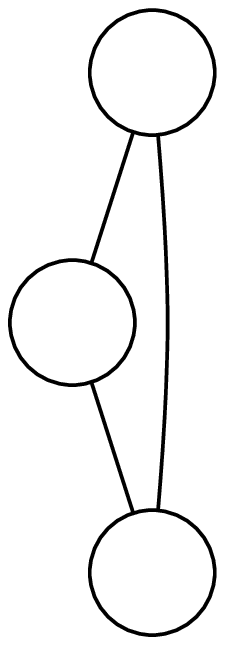}
\includegraphics[scale=0.3]{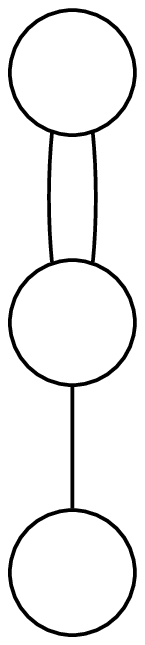}
\includegraphics[scale=0.3]{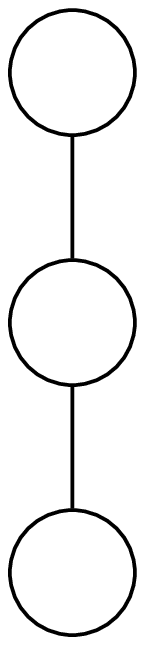}
\caption{The 3 C-trees with 3 nodes.}
\label{fig.R3}
\end{figure}

(i) The skeleton the tree on one node generates one C-tree
with a cycle of length 3,
$r_1(3)=1$.
(ii) The skeleton tree of two nodes generates one C-tree
(weights 1+2), $r_2(3)=1$.
(iii) The skeleton tree of three nodes generates one C-tree, $r_3(3)=1$
(Figure \ref{fig.R3}) :
\begin{equation}
c(3)=1+1+1=3.
\end{equation}
\subsection{Four nodes}
\begin{figure}
\includegraphics[scale=0.3]{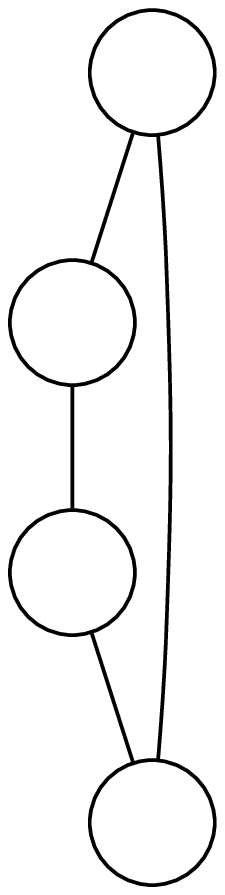}

\includegraphics[scale=0.3]{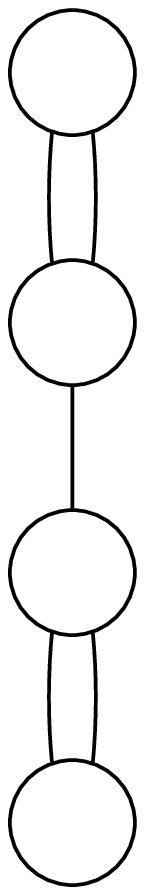}
\includegraphics[scale=0.3]{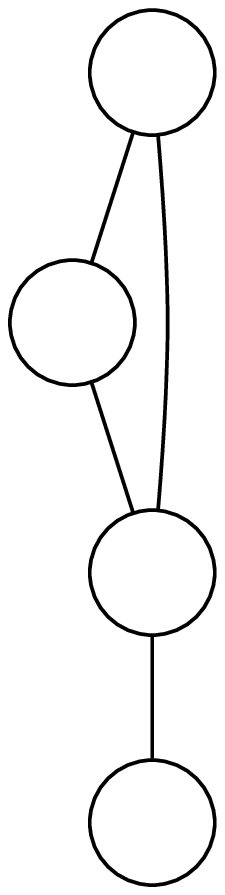}

\includegraphics[scale=0.3]{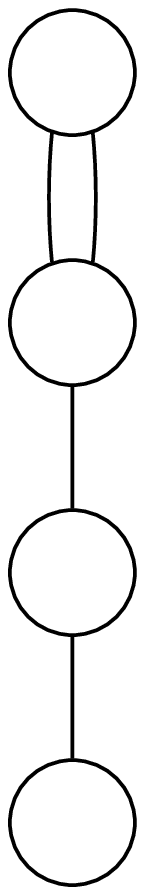}
\includegraphics[scale=0.3]{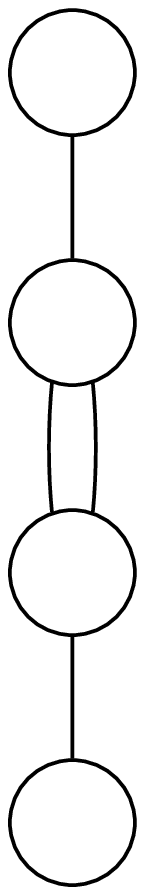}
\includegraphics[scale=0.3]{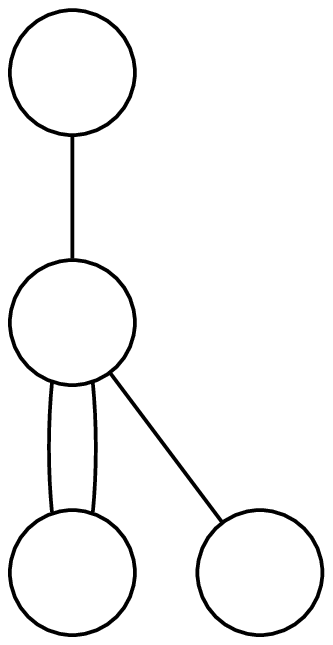}

\includegraphics[scale=0.3]{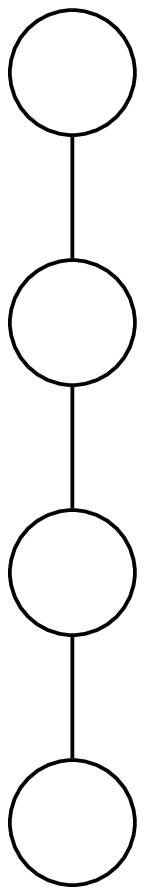}
\includegraphics[scale=0.3]{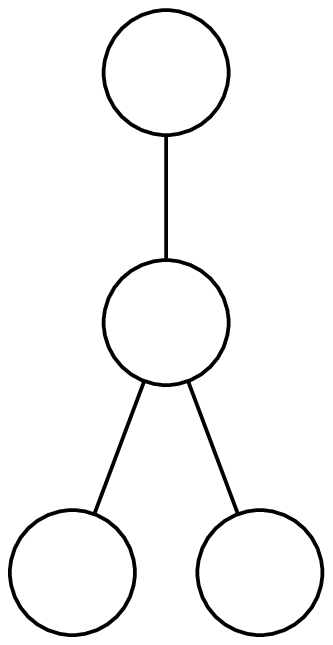}

\caption{The 8 C-trees with 4 nodes.}
\label{fig.R4}
\end{figure}

(i) The skeleton tree on one node generates a C-tree
with a cycle of length 4 (top row Fig.\ \ref{fig.R4}):
$r_1(4)=1$.
(ii) The skeleton tree of two nodes generates two C-trees
(weights 1+3=2+2, 2nd row Figure \ref{fig.R4}): $r_2(4)=2$
(iii) The skeleton tree of three nodes generates one C-tree
with a cycle at a leave (weights 2+1+1),
and two C-trees with cycle at the center node (weights 1+2+1), where the
leaves are either attached to the same cycle node or to two different
cycle nodes. This turns out to be the simplest example where mapping
a weighted skeleton tree generates more than one C-tree,
and where $r_t(n)$ surpasses the number of weighted trees of weight $n$
with $t$ nodes \cite[A303841]{EIS}.
Summarized in the 3rd row of Figure \ref{fig.R4}: $r_3(4)=3$.
(iv) The two skeleton trees of 4 nodes are C-trees as they are
(last row Figure \ref{fig.R4}): $r_4(4)=2$.
\begin{equation}
c(4)=\sum_{t=1}^4 r_t(4) = 1+2+3+2=8.
\end{equation}

\subsection{Five nodes}
\begin{figure}

\includegraphics[scale=0.3]{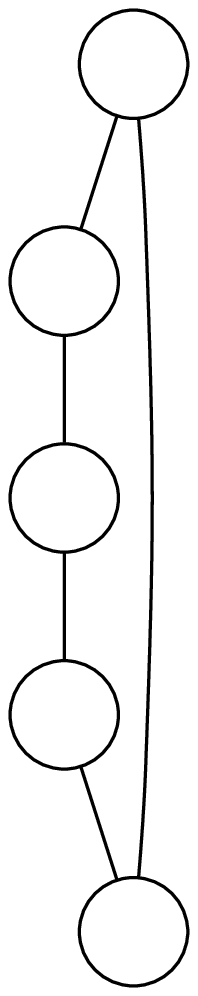}

\includegraphics[scale=0.3]{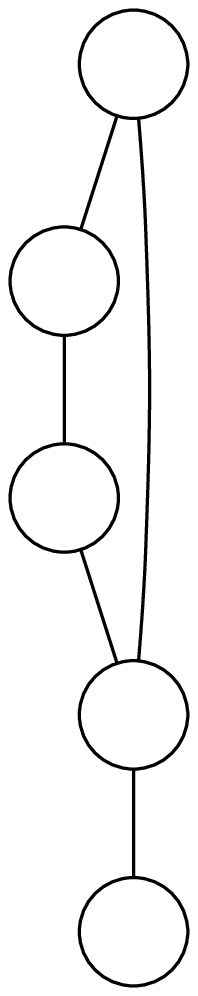}
\includegraphics[scale=0.3]{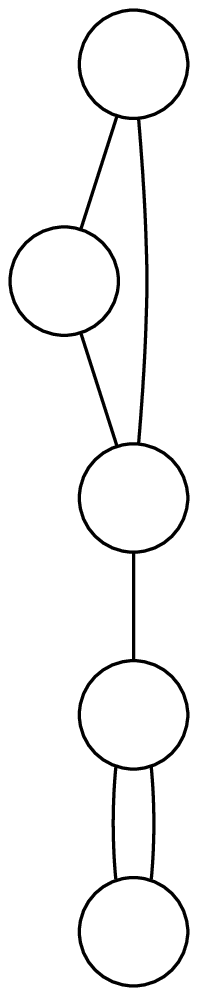}

\includegraphics[scale=0.3]{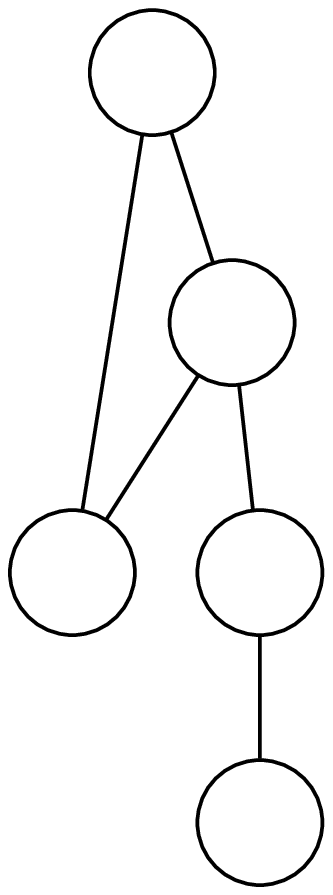}
\includegraphics[scale=0.3]{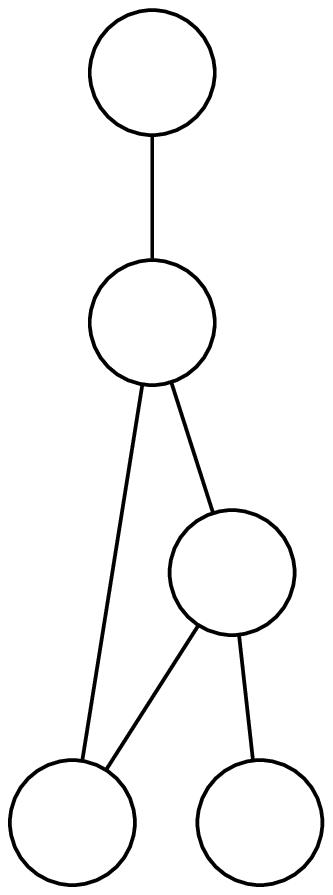}
\includegraphics[scale=0.3]{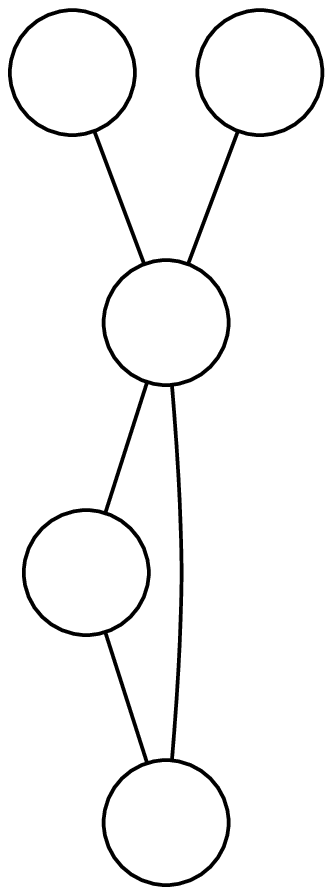}
\includegraphics[scale=0.3]{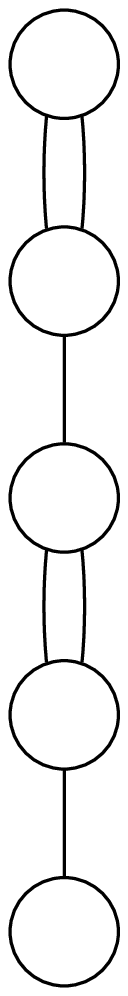}
\includegraphics[scale=0.3]{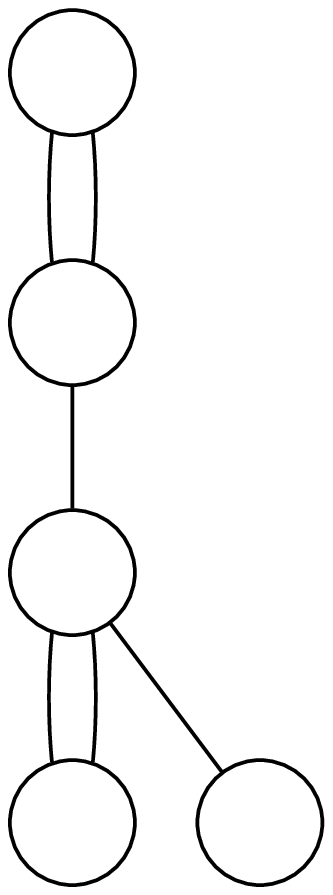}
\includegraphics[scale=0.3]{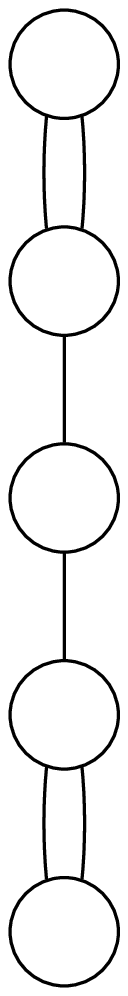}

\includegraphics[scale=0.3]{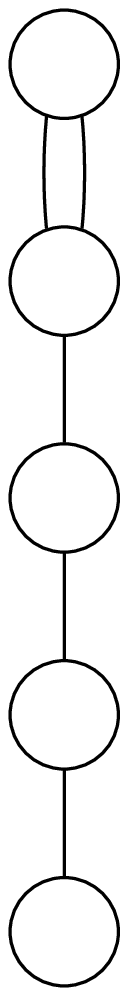}
\includegraphics[scale=0.3]{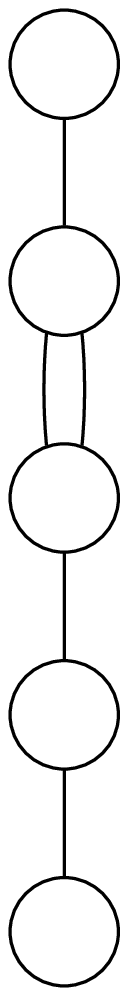}
\includegraphics[scale=0.3]{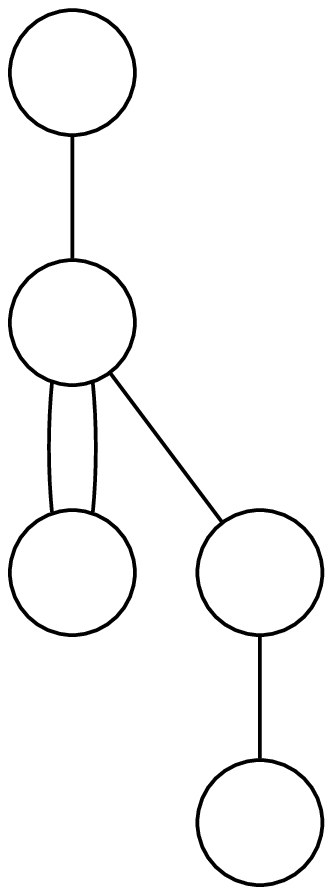}

\includegraphics[scale=0.3]{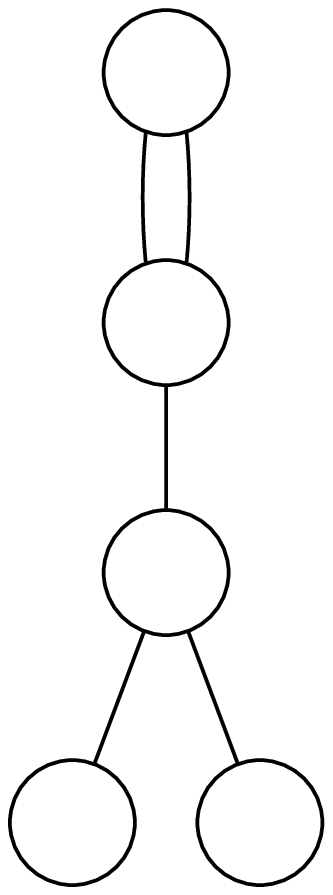}
\includegraphics[scale=0.3]{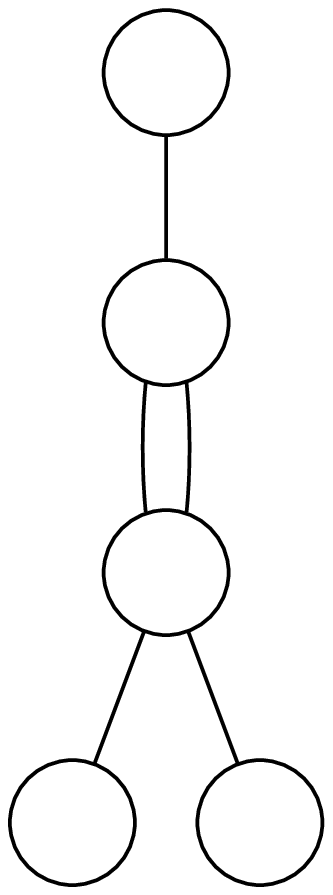}
\includegraphics[scale=0.3]{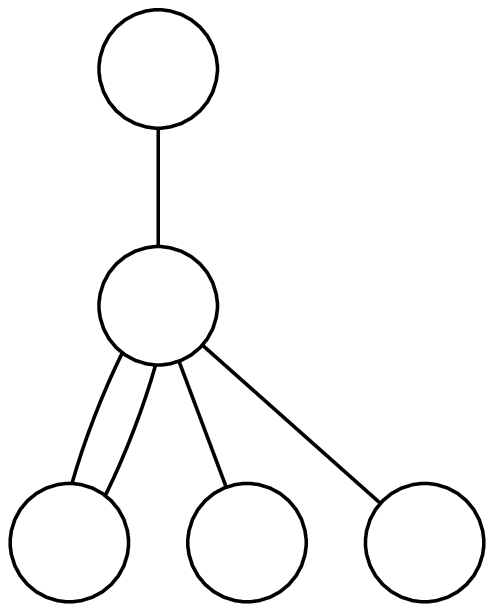}

\includegraphics[scale=0.3]{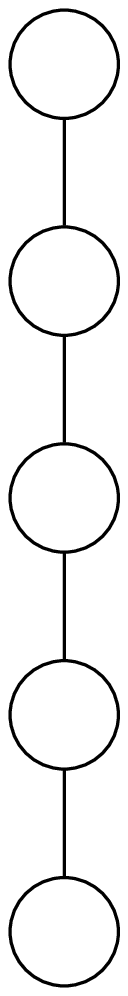}
\includegraphics[scale=0.3]{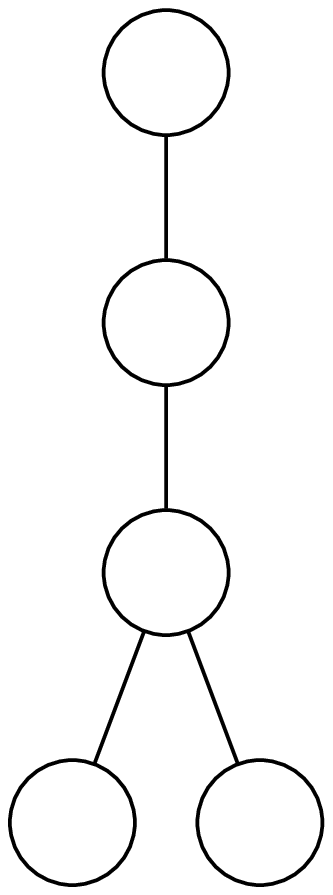}
\includegraphics[scale=0.3]{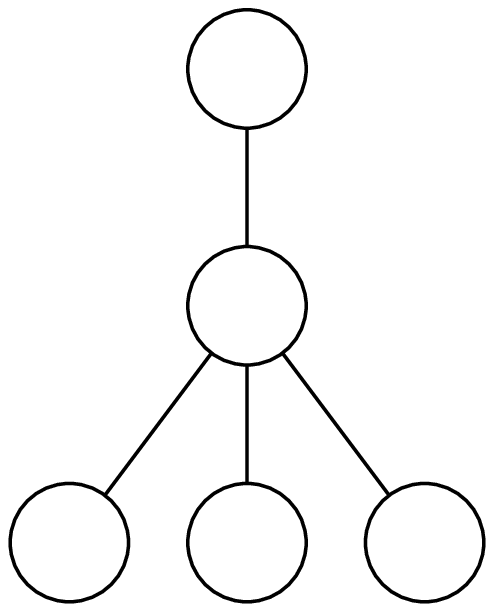}

\caption{The 18 C-trees with 5 nodes, weakly sorted by
the number of nodes in their skeleton trees.}
\label{fig.R5}
\end{figure}

(i) The skeleton tree with one node generates one C-tree with the cycle of 5
(top row Figure \ref{fig.R5}):
$r_1(5)=1$.
(ii) The skeleton tree with two nodes and weights $1+4=2+3$
generates (2nd row Figure \ref{fig.R5}) $r_2(5)=2$.
(iii) The skeleton tree with 3 nodes with weights $3+1+1=1+3+1=2+2+1=2+1+2$
generates a cycle of 3 on one leaf, 
a cycle of 3 in the middle,
cycles of 2 at one leaf and in the middle,
or two cycles of 2 on both leaves (3rd row Figure \ref{fig.R5}):
$r_3(5)=1+2+2+1=6$.
(iv) The linear skeleton tree with 4 nodes with weights $2+1+1+1=1+2+1+1$
generates the cycle either at a leaf or at one of the bi-centers
(4th row Figure \ref{fig.R5}): $r_{4a}(5)=1+2=3$.
The star skeleton tree with 4 nodes generates the cycle of 2 either at a leaf
(1 C-tree)
or in the middle (2 C-trees) (5th row Figure \ref{fig.R5}): $r_{4b}(5)=1+2=3$.
$r_4(5)=r_{4a}(5)+r_{4b}(5)=6$.
(v) The skeleton trees with 5 nodes represent C-trees as they are
(last row Figure \ref{fig.R5}): $r_5(5)=3$.
\begin{equation}
c(5)= \sum_{t=1}^5 r_t(5) = 1+2+6+6+3=18.
\end{equation}
\clearpage
\subsection{Six nodes}
\begin{figure}

\includegraphics[scale=0.3]{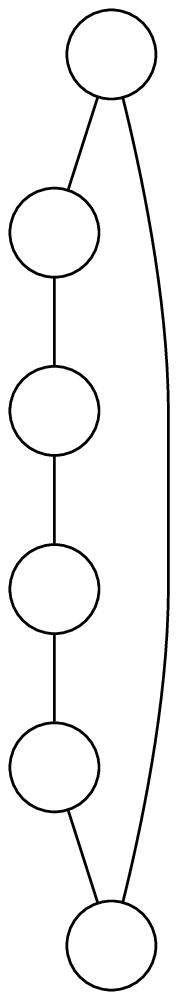}

\includegraphics[scale=0.3]{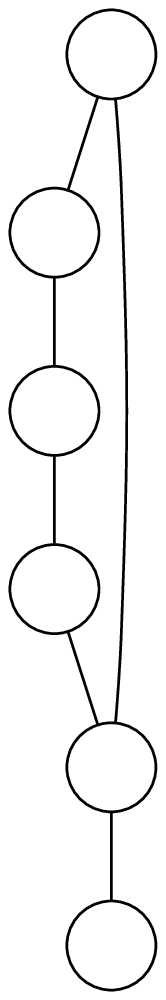}
\includegraphics[scale=0.3]{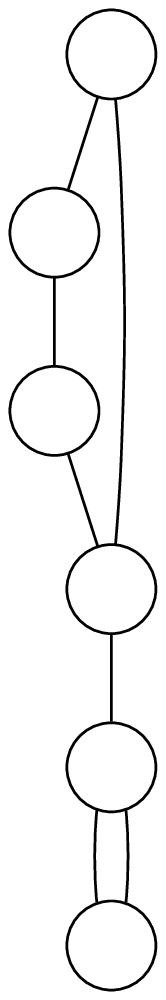}
\includegraphics[scale=0.3]{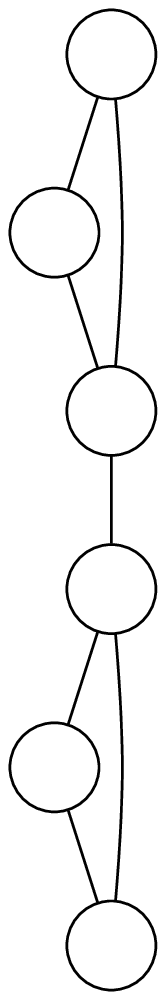}

\includegraphics[scale=0.3]{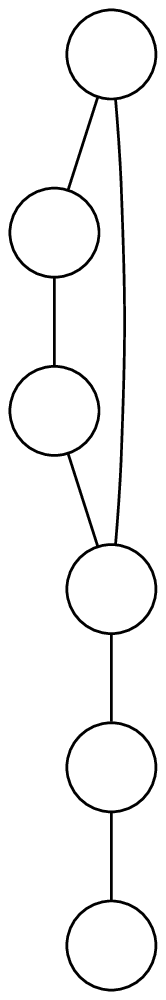}
\includegraphics[scale=0.3]{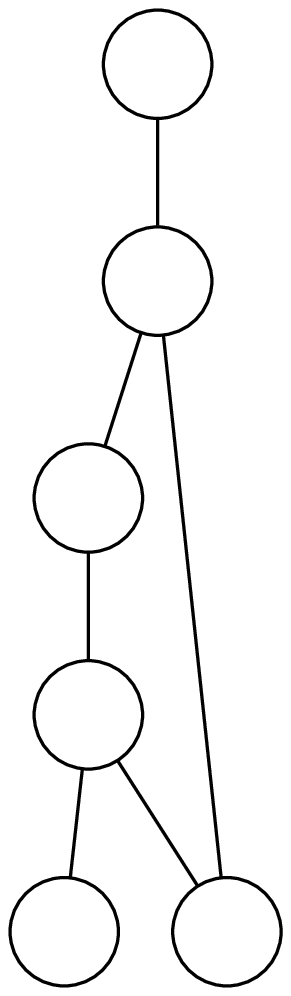}
\includegraphics[scale=0.3]{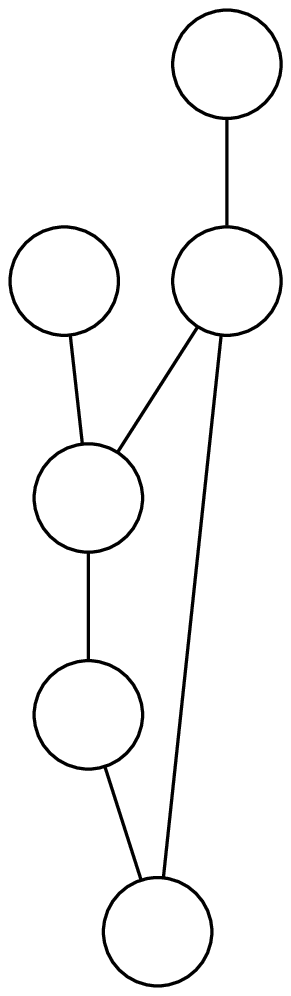}
\includegraphics[scale=0.3]{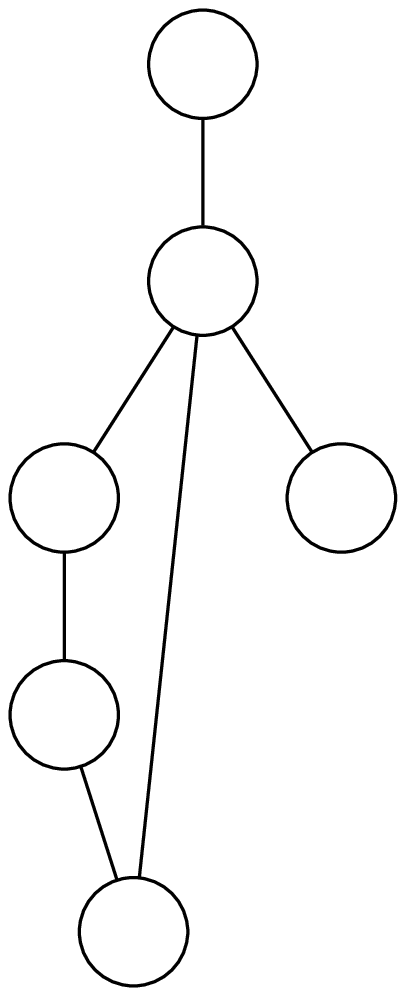}
\includegraphics[scale=0.3]{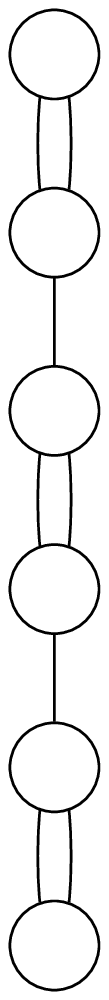}
\includegraphics[scale=0.3]{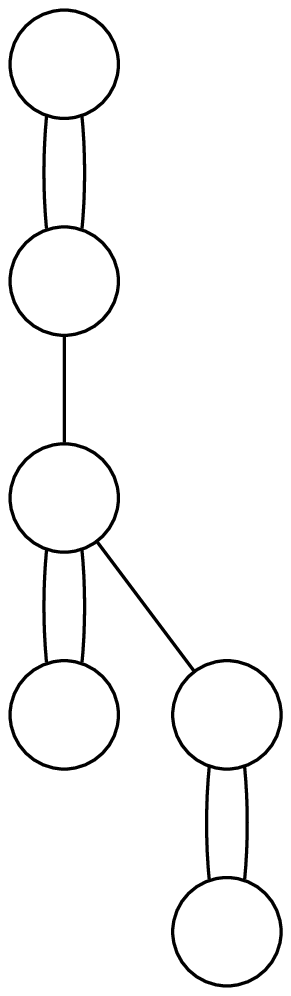}
\includegraphics[scale=0.3]{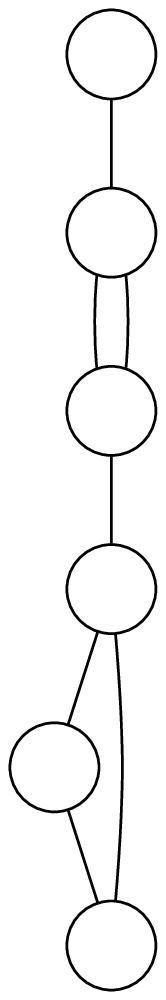}
\includegraphics[scale=0.3]{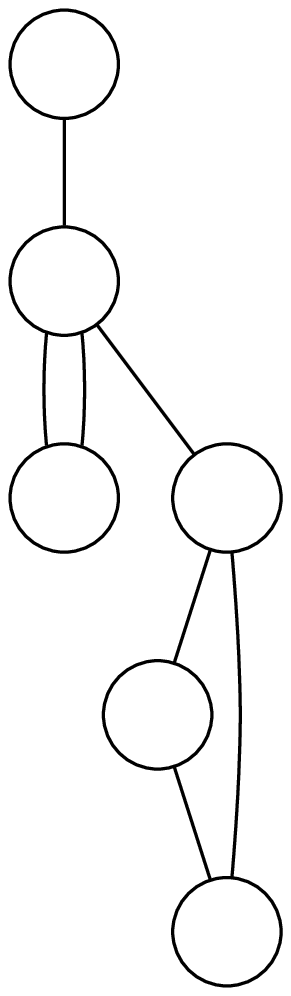}
\includegraphics[scale=0.3]{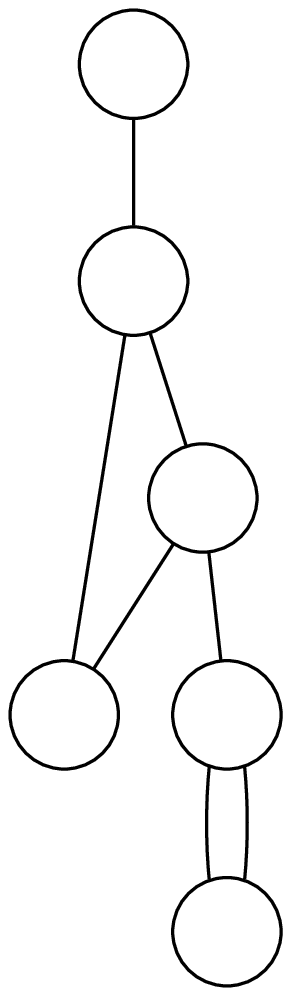}
\includegraphics[scale=0.3]{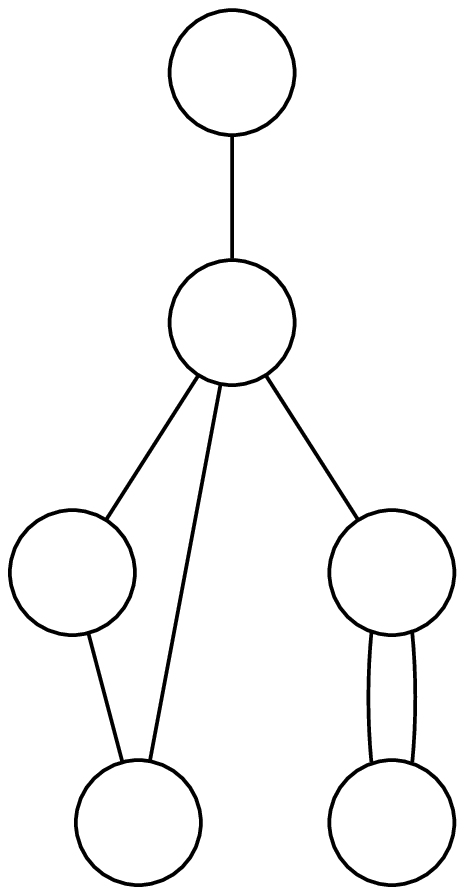}
\includegraphics[scale=0.3]{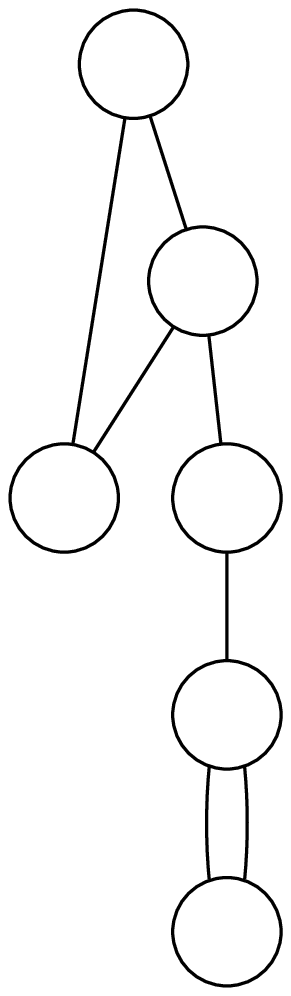}
\caption{The $r_1(6)=1$, $r_2(6)=3$ and $r_3(6)=11$ C-trees with 6 nodes and skeleton trees of 1, 2  or 3 nodes.}
\label{fig.R6a}
\end{figure}

(i) The skeleton tree with with node generates one C-tree with a cycle of 6
(top row Figure \ref{fig.R6a}):
$r_1(6)=1$.

(ii) The skeleton tree with 2 nodes 
generates
$r_2(6)=3$ C-trees (2nd row Figure \ref{fig.R6a}).

(iii) For the skeleton tree with 3 nodes consider weights
$4+1+1=1+4+1=2+2+2=1+2+3=1+3+2=3+1+2$. 
They generate in that order
$r_3(6)=1+3+2+2+2+1=11$ C-trees (last row Figure \ref{fig.R6a}).

\begin{figure}

\includegraphics[scale=0.3]{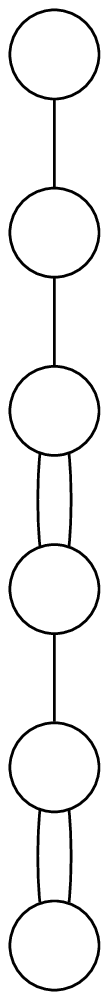}
\includegraphics[scale=0.3]{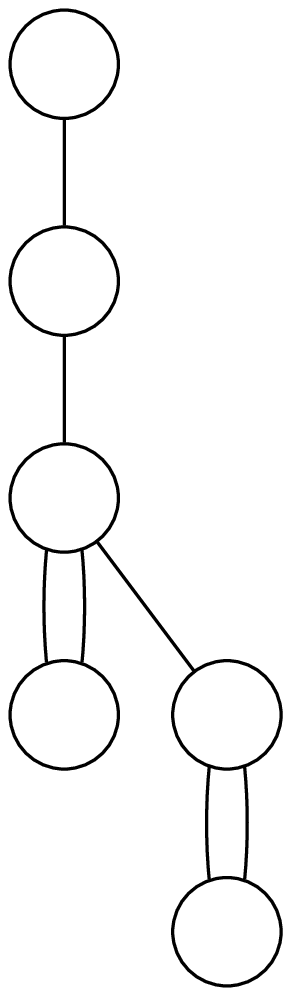}
\includegraphics[scale=0.3]{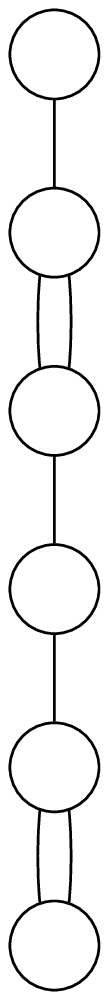}
\includegraphics[scale=0.3]{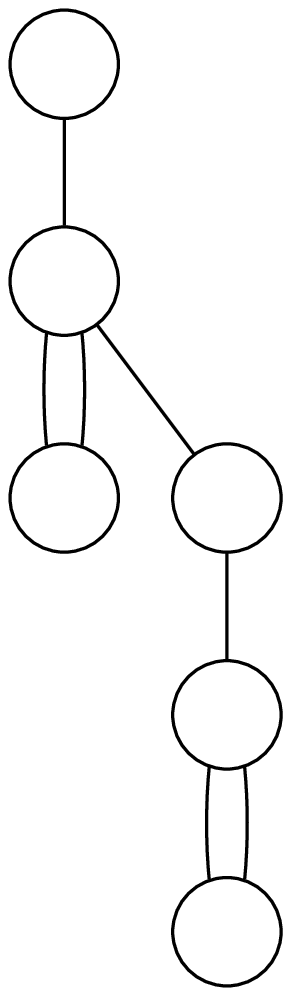}
\includegraphics[scale=0.3]{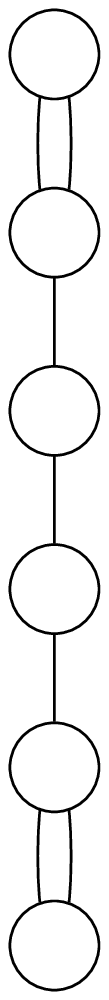}
\includegraphics[scale=0.3]{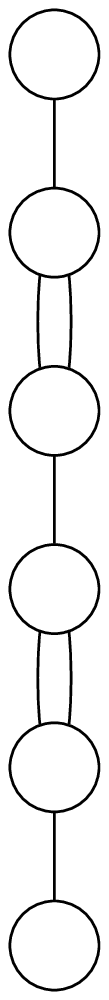}
\includegraphics[scale=0.3]{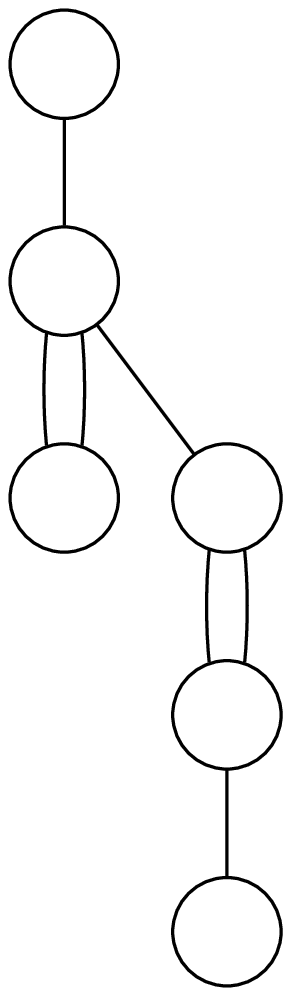}
\includegraphics[scale=0.3]{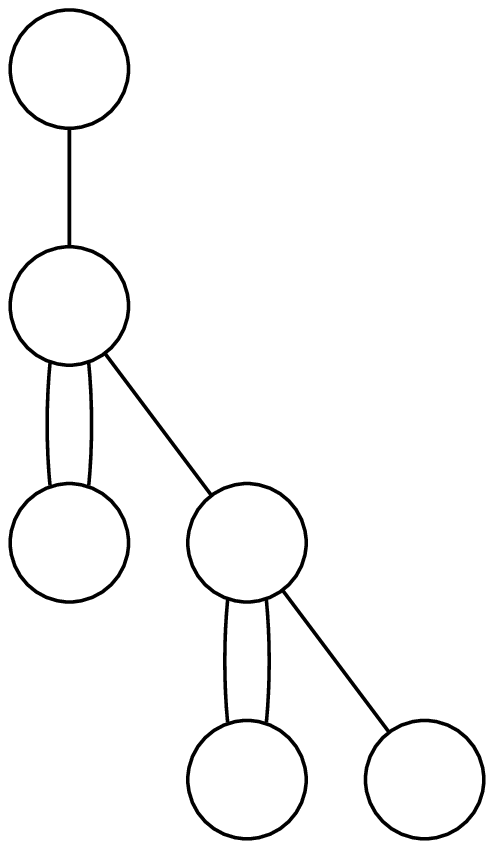}
\includegraphics[scale=0.3]{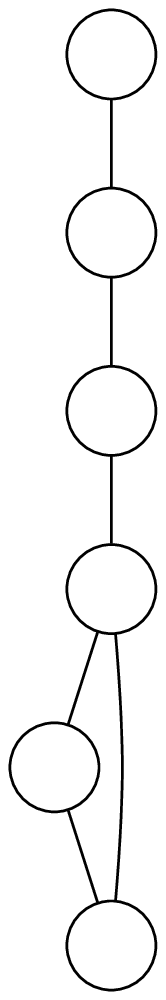}
\includegraphics[scale=0.3]{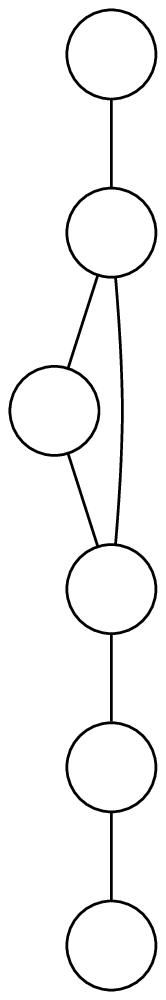}
\includegraphics[scale=0.3]{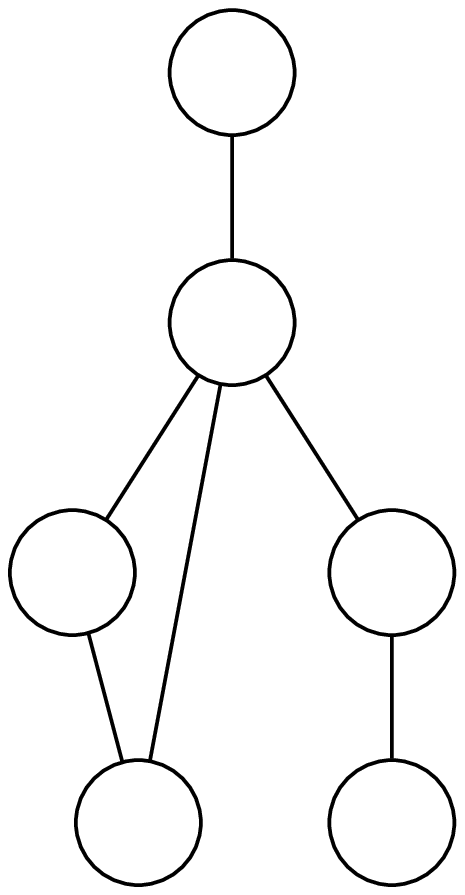}

\includegraphics[scale=0.3]{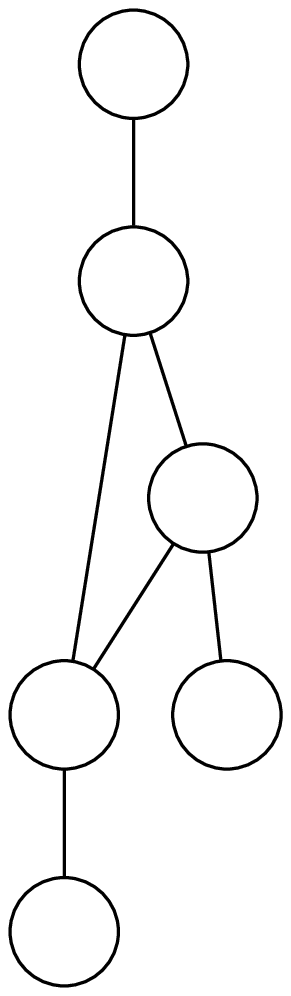}
\includegraphics[scale=0.3]{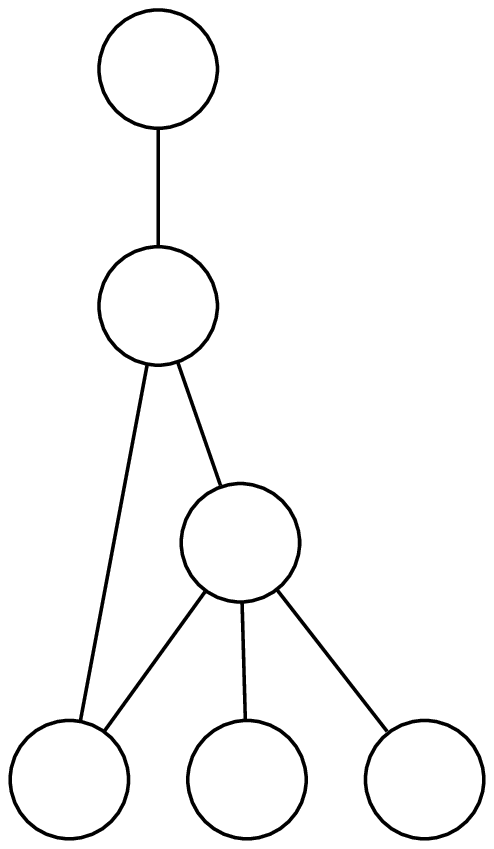}
\includegraphics[scale=0.3]{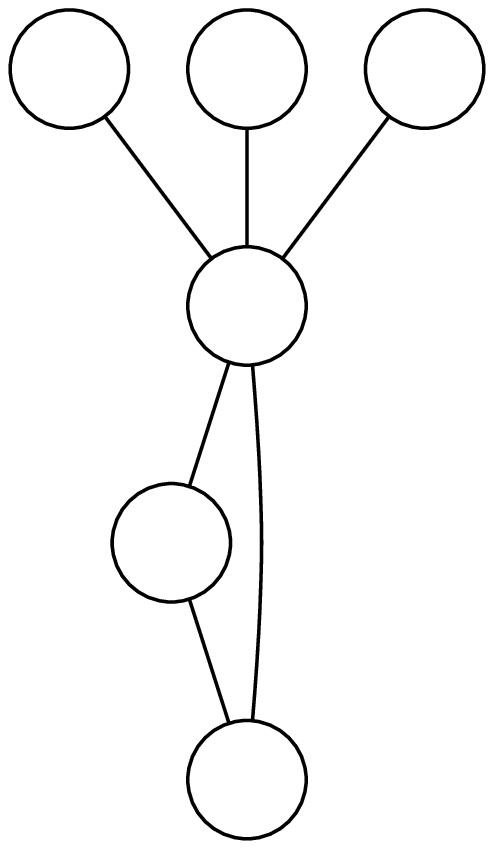}
\includegraphics[scale=0.3]{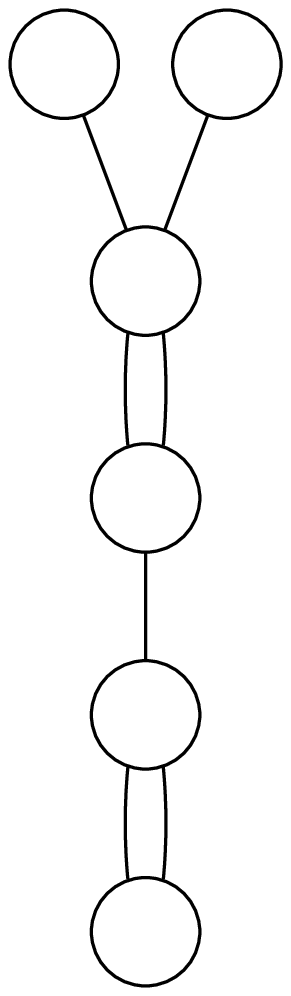}
\includegraphics[scale=0.3]{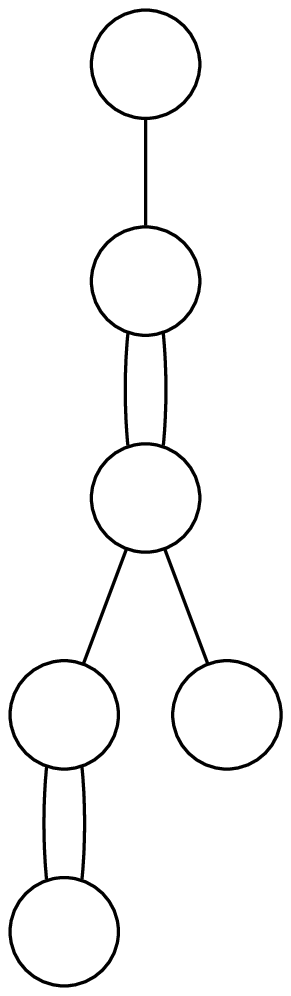}
\includegraphics[scale=0.3]{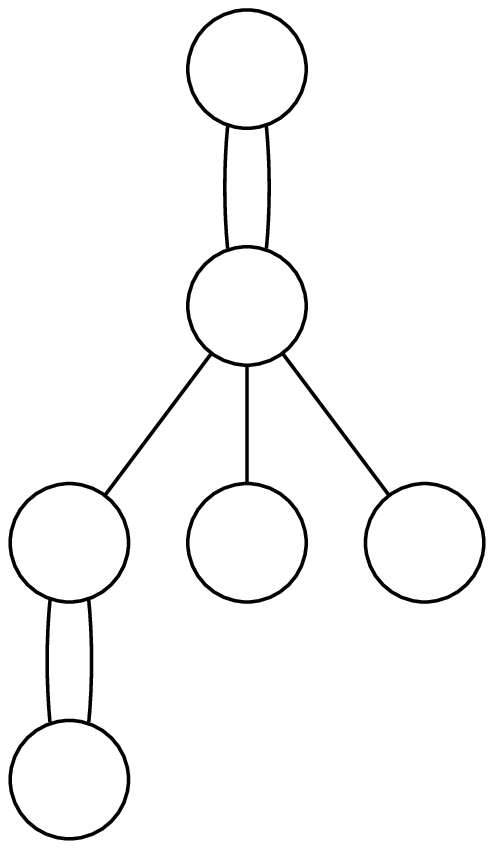}
\includegraphics[scale=0.3]{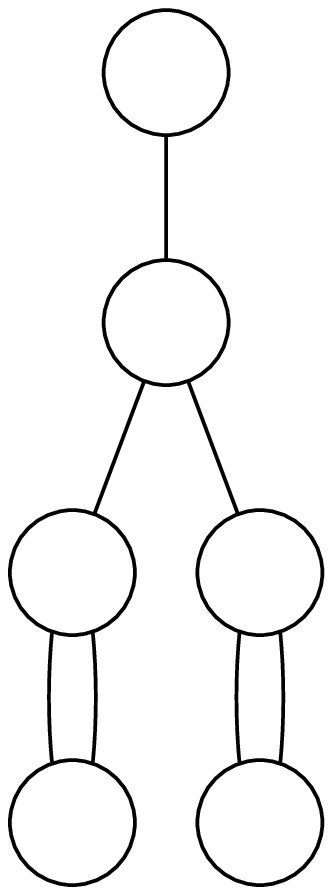}
\includegraphics[scale=0.3]{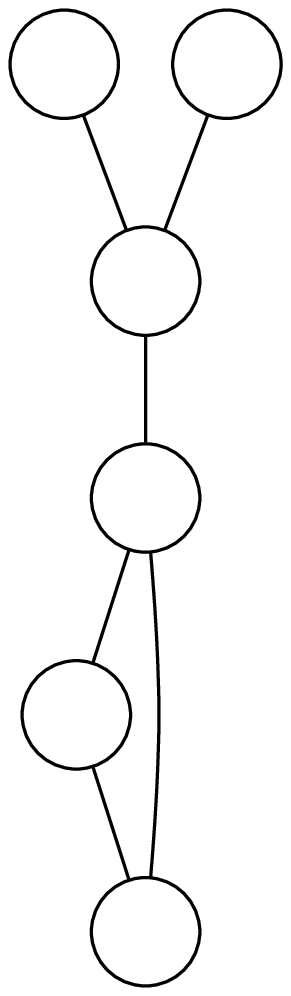}

\includegraphics[scale=0.3]{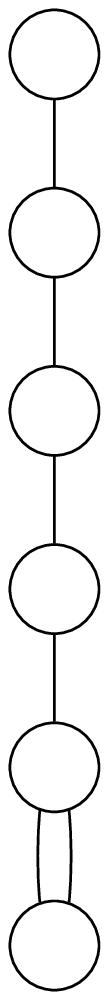}
\includegraphics[scale=0.3]{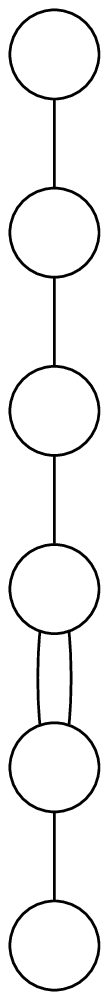}
\includegraphics[scale=0.3]{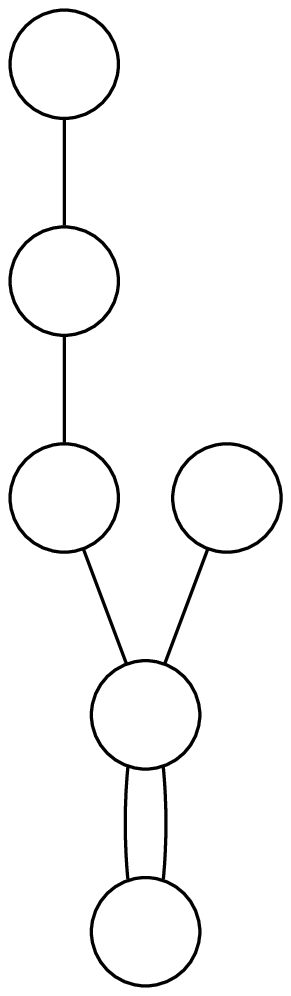}
\includegraphics[scale=0.3]{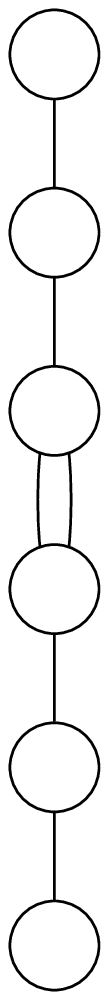}
\includegraphics[scale=0.3]{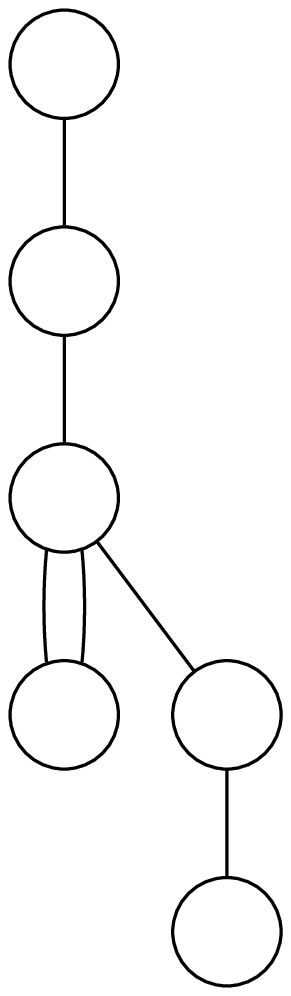}

\includegraphics[scale=0.3]{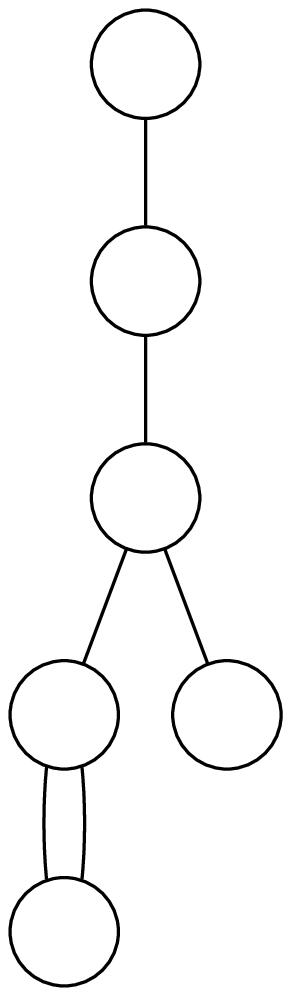}
\includegraphics[scale=0.3]{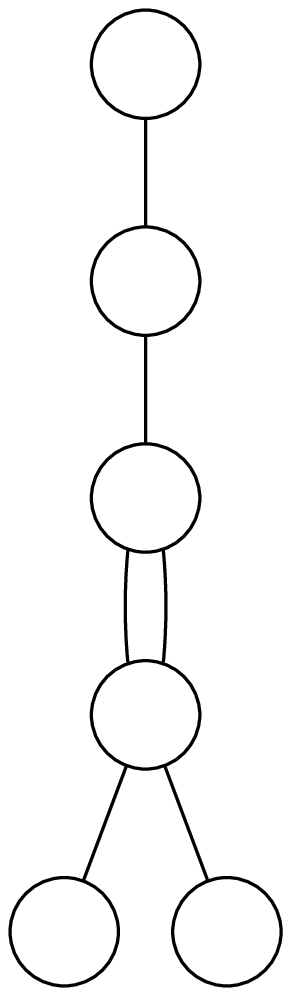}
\includegraphics[scale=0.3]{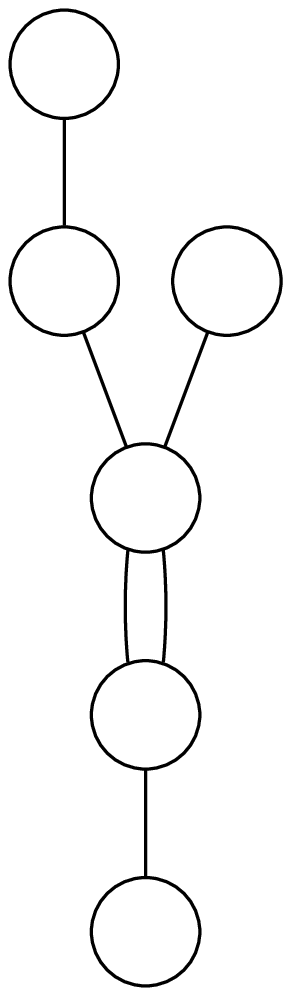}
\includegraphics[scale=0.3]{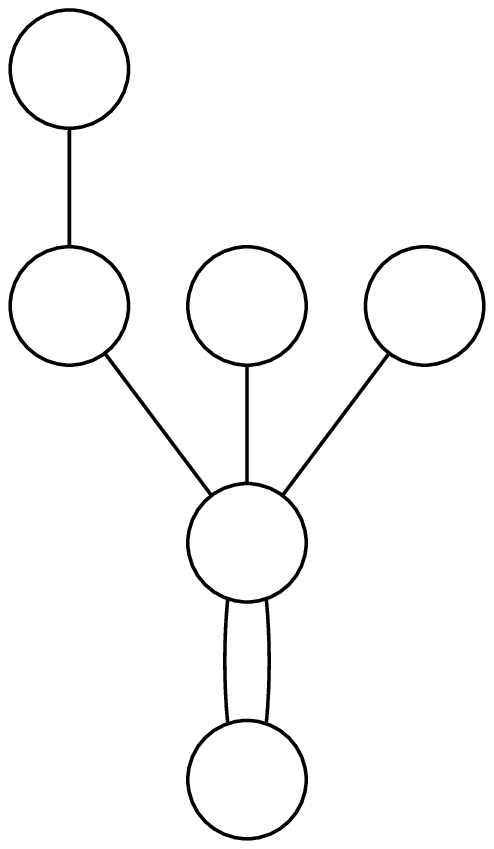}
\includegraphics[scale=0.3]{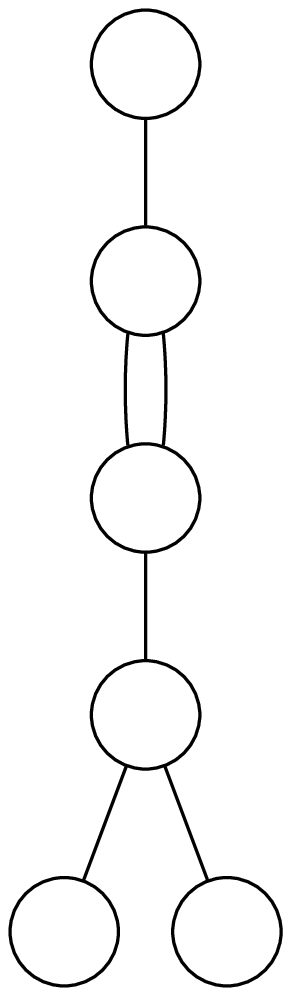}
\includegraphics[scale=0.3]{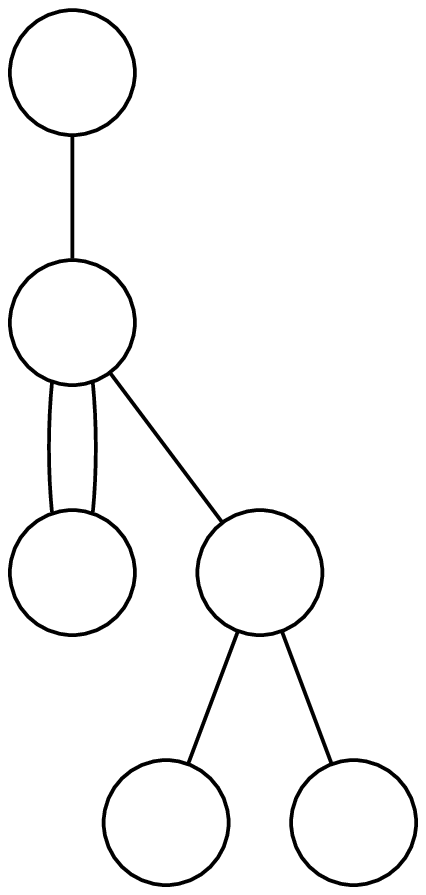}
\includegraphics[scale=0.3]{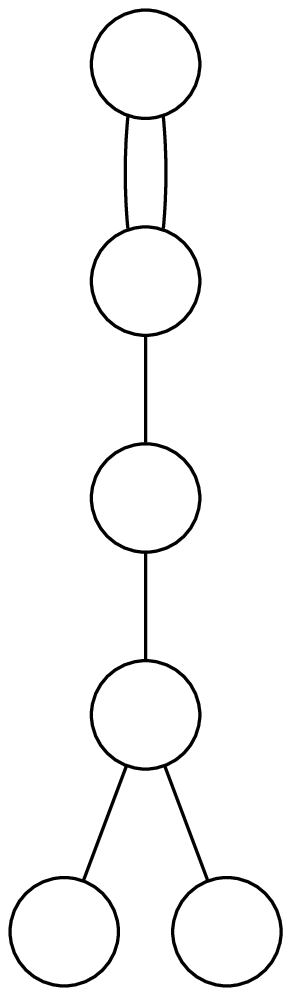}

\includegraphics[scale=0.3]{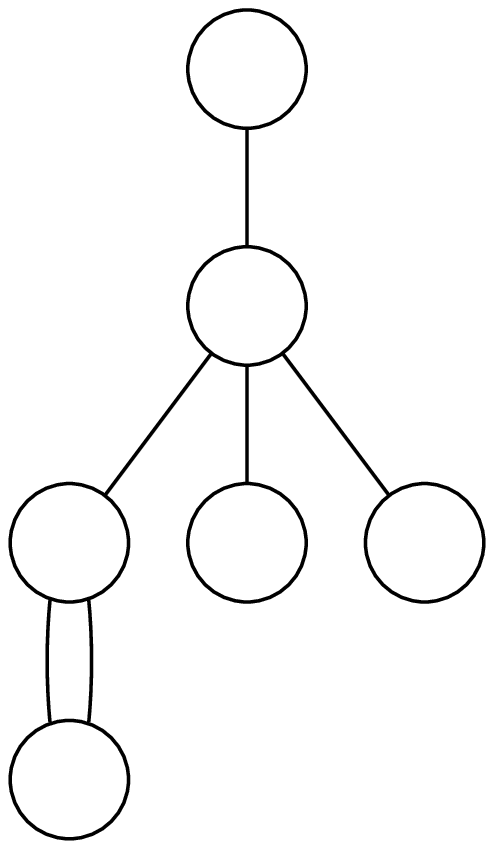}
\includegraphics[scale=0.3]{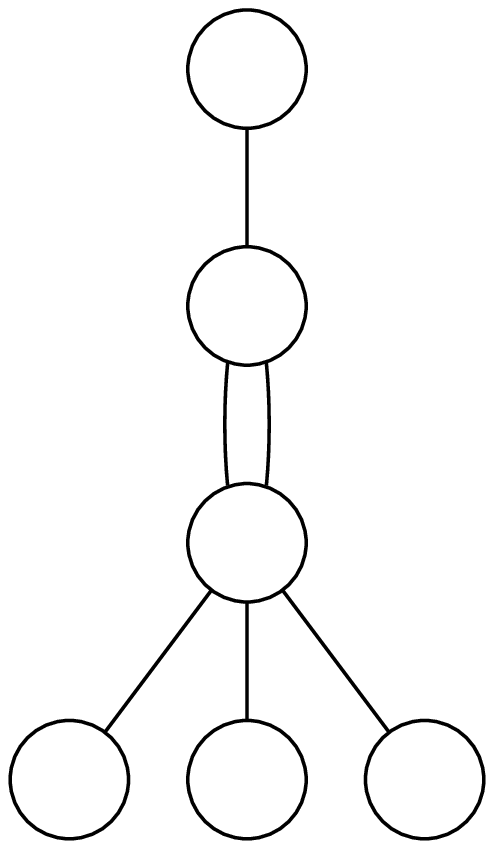}
\includegraphics[scale=0.3]{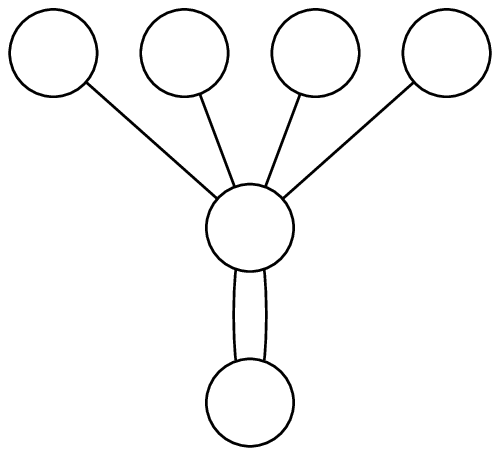}
\includegraphics[scale=0.3]{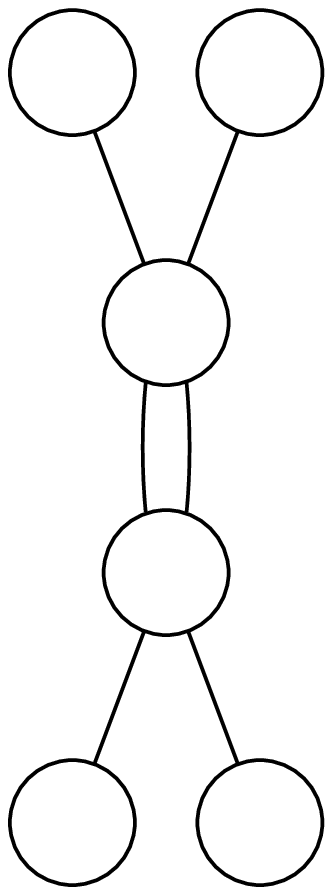}

\includegraphics[scale=0.3]{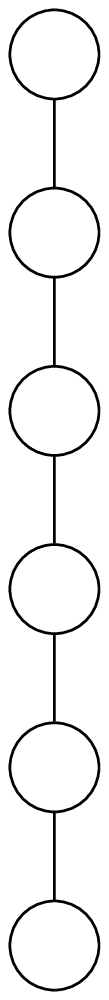}
\includegraphics[scale=0.3]{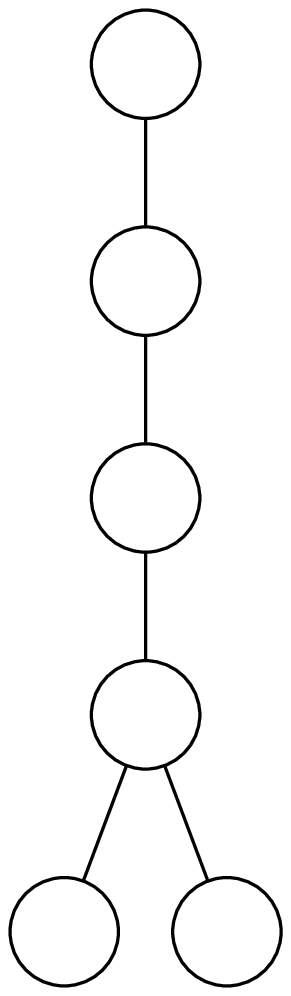}
\includegraphics[scale=0.3]{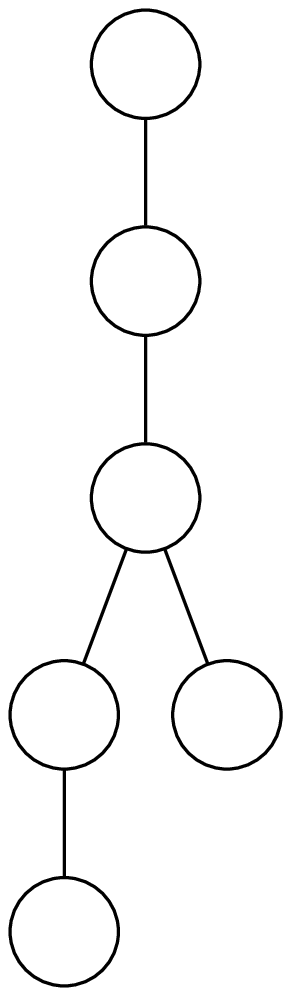}
\includegraphics[scale=0.3]{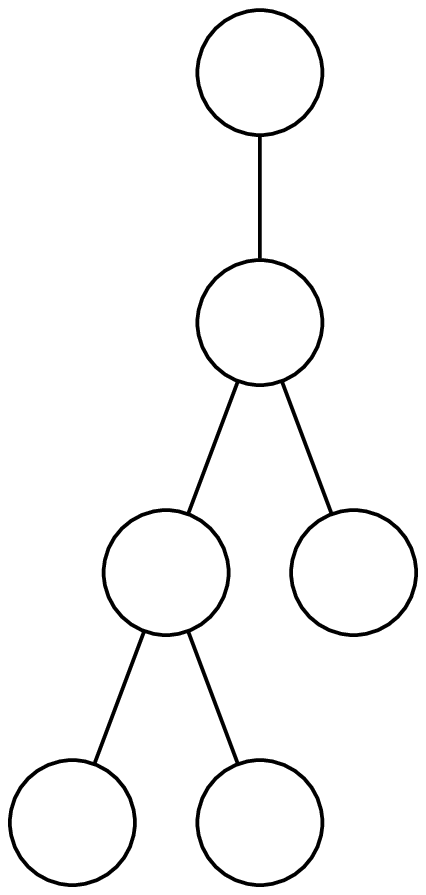}
\includegraphics[scale=0.3]{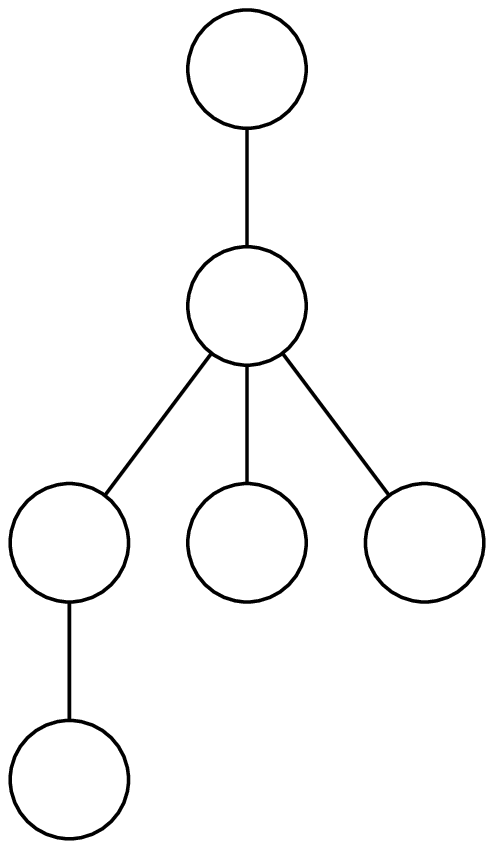}
\includegraphics[scale=0.3]{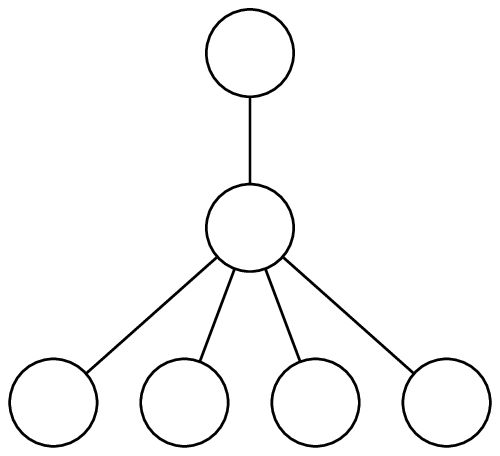}
\caption{The $r_4(6)=19$, $r_5(6)=16$ and $r_6(6)=6$ C-trees with 6 nodes (skeleton trees of 4, 5 or 6 nodes).}
\label{fig.R6b}
\end{figure}

(iv) For the linear skeleton tree with 4 nodes and weight 6 consider the
weights $1+1+2+2 = 1+2+1+2 = 2+1+1+2 = 1+2+2+1 = 1+1+1+3 = 1+1+3+1$
for the nodes left-to-right,
meaning two cycles of 2 are adjacent or not, at leaves or mixed, and
the cycle of 3 is a leaf or a center.
Illustrated in that order by the top row of Figure \ref{fig.R6b} these generate
$r_{4a}(6)=2+2+1+3+1+2=11$ C-trees.

For the star skeleton tree with 4 nodes and weight 6 consider the weight
partitions
$3+1+1+1 = 2+1+1+2 = 1+2+2+1=1+3+1+1$, where the 3-cycle, 2-cycle or 1-cycle 
is in the center. 
Illustrated in that order in the second row of Figure \ref{fig.R6b} this
generates
$r_{4b}(6)=3+3+1+1=8$ C-trees.
$r_4(6)=r_{4a}(6)+r_{4b}(6)=11+8=19$.

(v) From the linear skeleton tree with 5 nodes and weight 6 consider
the partitions $1+1+1+1+2=1+1+1+2+1 =1+1+2+1+1$ with a 2-cycle at one
of the 3 positions. $1+1+1+1+2$ generates one C-tree. $1+1+1+2+1$ generates
2
C-trees. $1+1+2+1+1$ generates 2 C-trees (third row of figure \ref{fig.R6b}):
$r_{5a}(6)=1+2+2=5$.

For the split skeleton tree with 5 nodes and weight 6 place the 2-cycle
either at one of the short ends (1 choice), at the node with degree 3
with 2 subtrees of 1 node and one subtree of 2 nodes (3 choices), at
the node with degree 2 (2 choices), or at the long end (1 choice)
(4th row Figure \ref{fig.R6b}):
$r_{5b}(6)=1+3+2+1=7$.

For the star skeleton tree with 5 nodes and weight 6 the 2-cycle can be at
a leaf (1 choice) or in the center (3 choices) according
to the 5th row of Figure \ref{fig.R6b}:
$r_{5c}(6)=1+3=4$.
The three types of skeleton trees with 5 nodes map to
$r_5(6)=r_{5a}(6)+r_{5b}(6)+r_{5c}(6)=5+7+4=16$ C-trees.

(vi) Each of the skeleton trees with 6 nodes generates one C-tree, $r_6(6)=6$
(last row Figure \ref{fig.R6b}).
\begin{equation}
c(6)=1+3+11+19+16+6=56.
\end{equation}

\begin{equation}
c(n) = 1,1,2,3,8,18,56,\ldots ;\quad n\ge 0
\label{eq.Rn}
\end{equation}
\clearpage

\section{Planted C-trees}
\begin{defn}
A planted C-tree is a C-tree with a marked node (called the root) which
is an endnode.
\end{defn}

The planted C-trees with $n$ nodes can be generated from the
C-trees with $n$ nodes by discarding the C-trees withou
endnodes, marking successively the endnodes of the remaining C-trees,
and retaining only these rooted C-trees which are nonequivalent
under the automorphism group of the graph.

\begin{defn}
The number of planted C-trees with $n$ nodes is $p(n)$,
and
\begin{equation}
P(x) =\sum_{n\ge 1 } p(n)x^n = x+x^2+2x^3+6x^4+\cdots
\end{equation}
the generating function. 
\end{defn}

\begin{figure}
\includegraphics[scale=0.3]{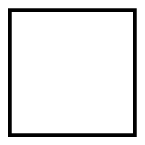}

\includegraphics[scale=0.3]{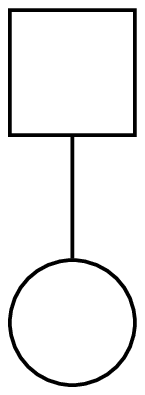}

\includegraphics[scale=0.3]{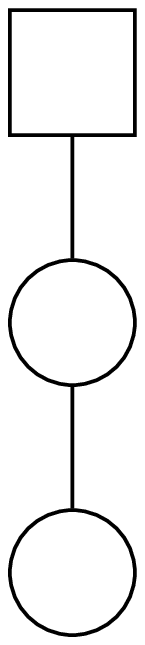}
\includegraphics[scale=0.3]{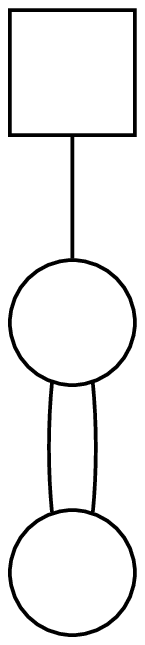}

\includegraphics[scale=0.3]{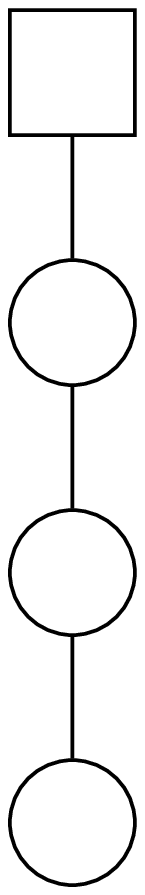}
\includegraphics[scale=0.3]{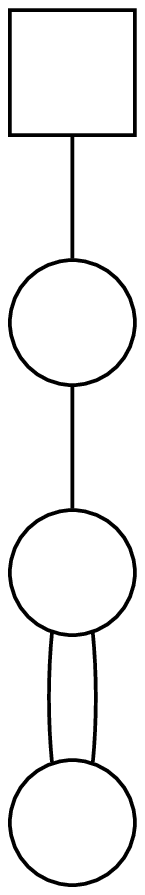}
\includegraphics[scale=0.3]{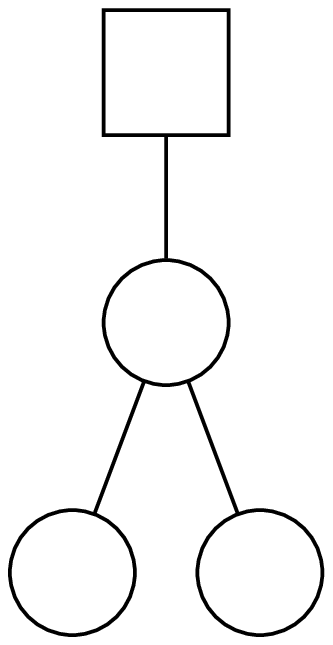}
\includegraphics[scale=0.3]{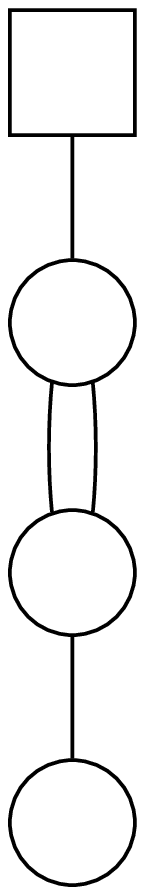}
\includegraphics[scale=0.3]{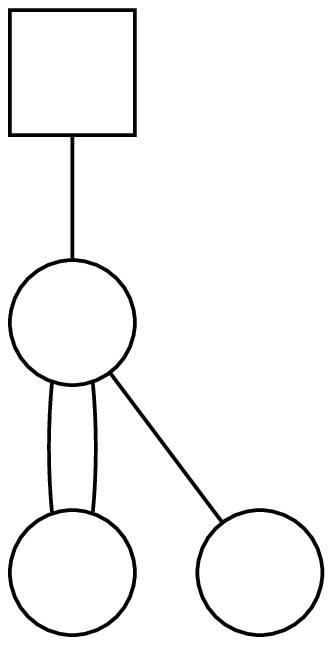}
\includegraphics[scale=0.3]{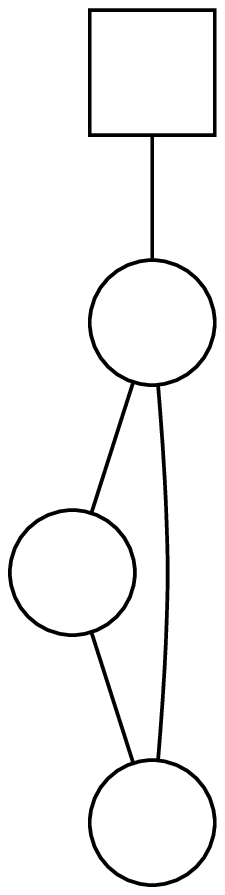}

\includegraphics[scale=0.3]{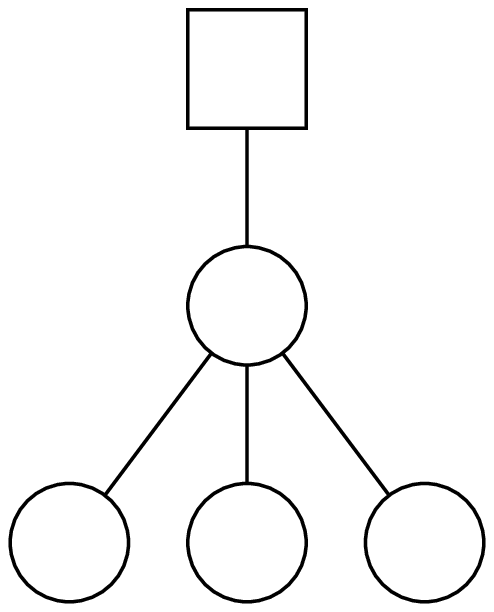}
\includegraphics[scale=0.3]{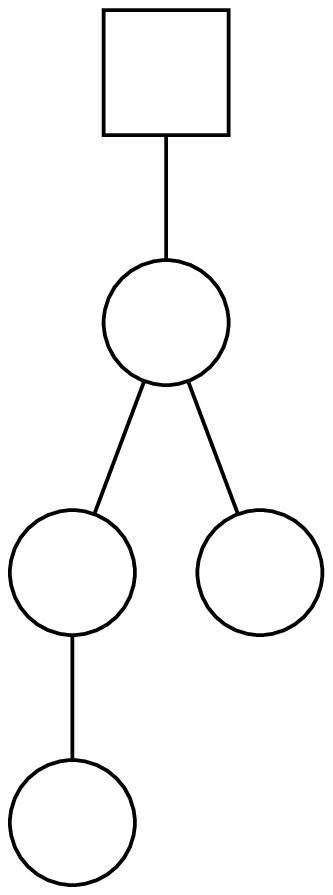}
\includegraphics[scale=0.3]{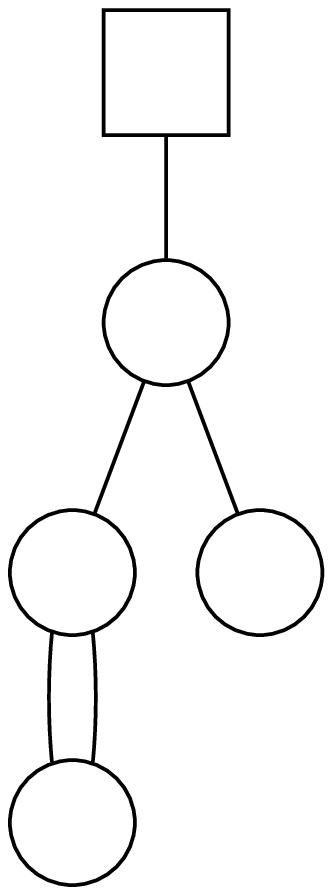}
\includegraphics[scale=0.3]{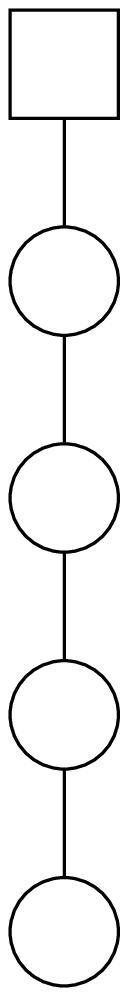}
\includegraphics[scale=0.3]{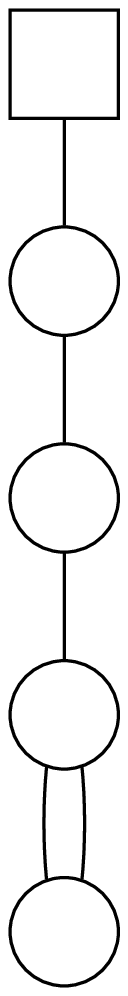}
\includegraphics[scale=0.3]{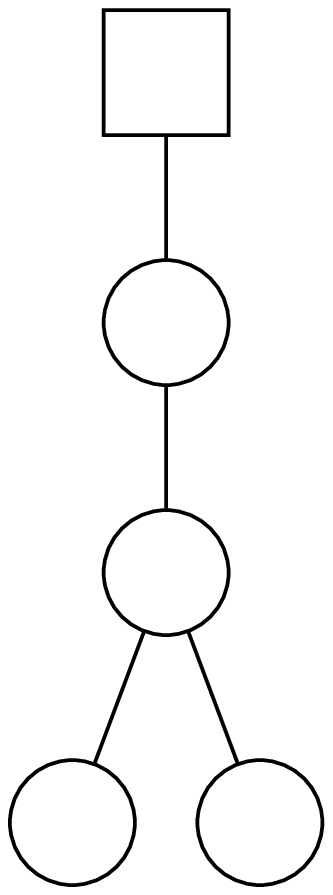}
\includegraphics[scale=0.3]{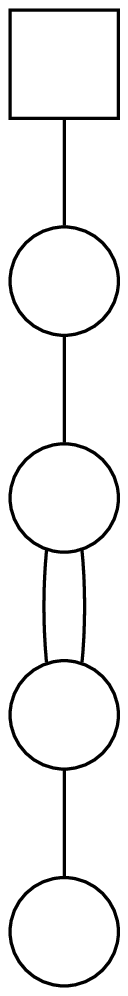}
\includegraphics[scale=0.3]{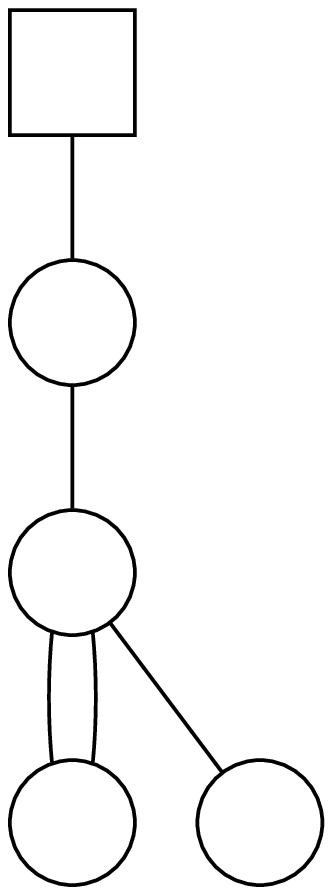}
\includegraphics[scale=0.3]{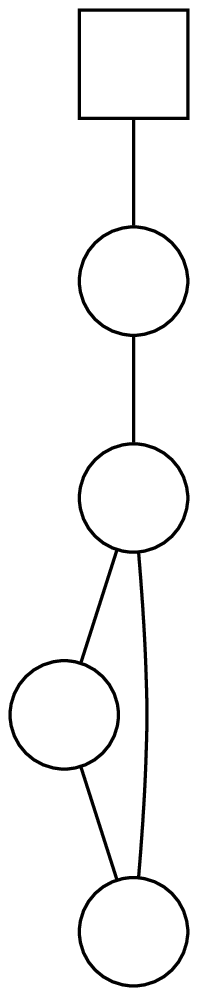}
\includegraphics[scale=0.3]{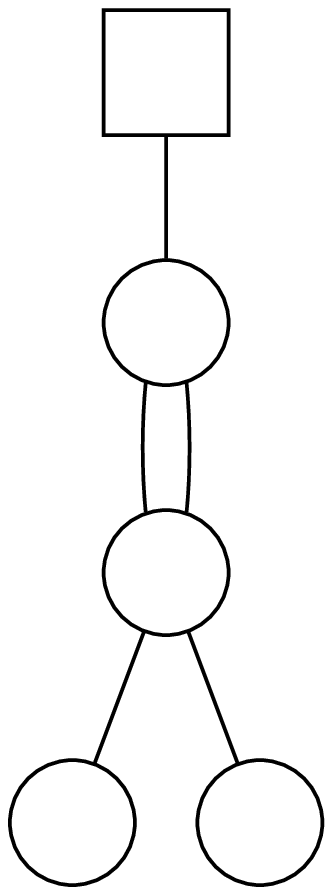}
\includegraphics[scale=0.3]{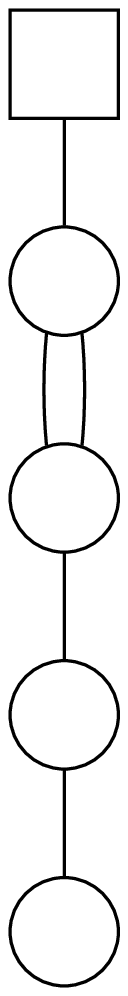}
\includegraphics[scale=0.3]{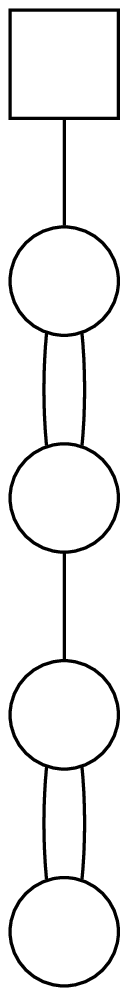}
\includegraphics[scale=0.3]{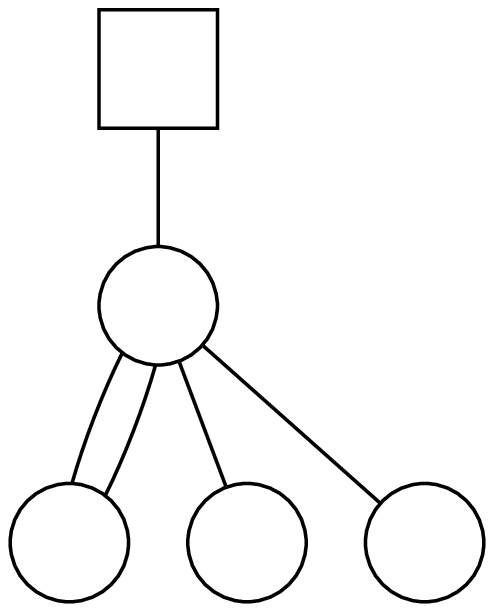}
\includegraphics[scale=0.3]{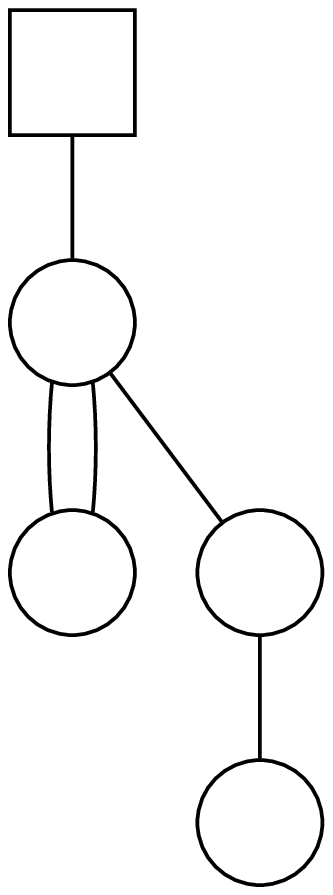}
\includegraphics[scale=0.3]{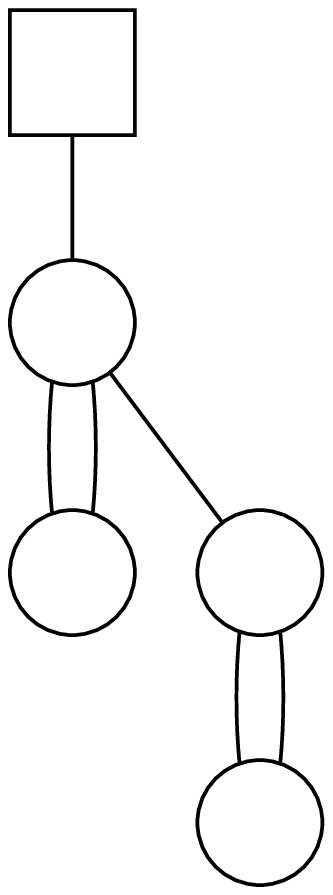}
\includegraphics[scale=0.3]{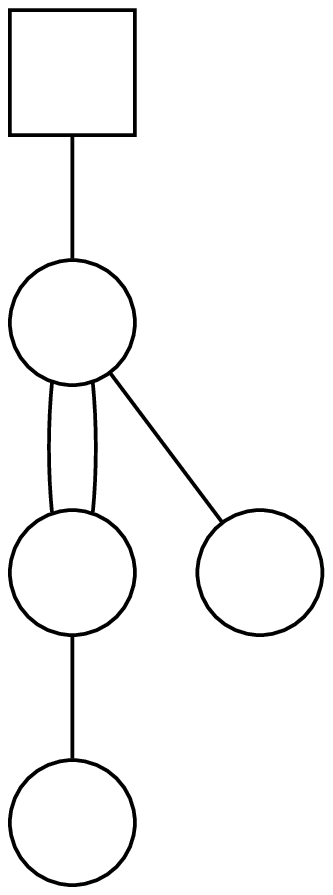}
\includegraphics[scale=0.3]{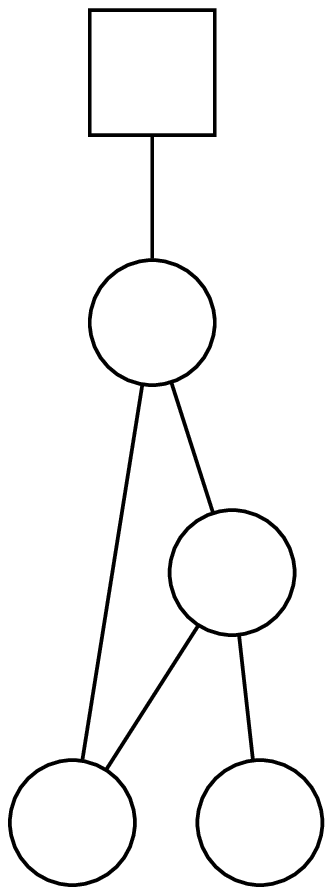}
\includegraphics[scale=0.3]{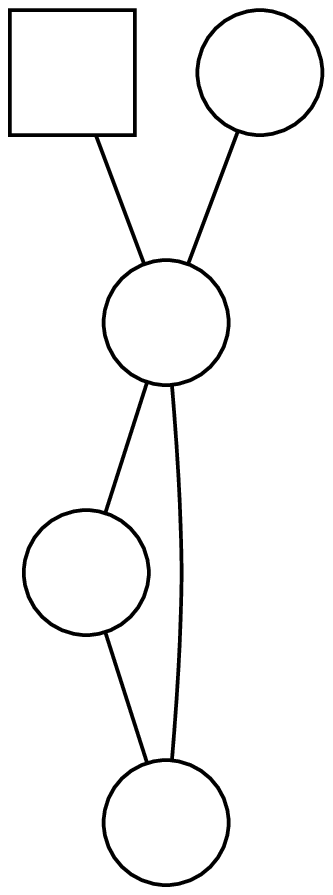}
\includegraphics[scale=0.3]{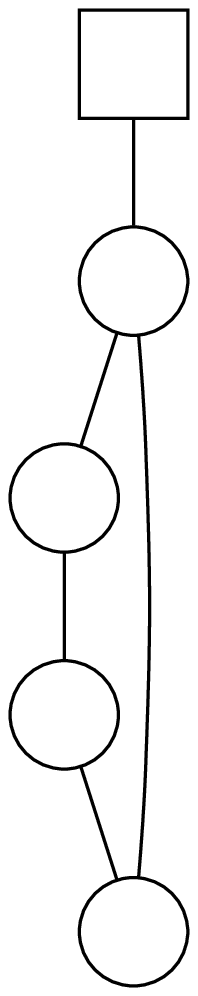}

\caption{Illustration of the planted C-trees up to 5 nodes, 1,1,2,6,19,\ldots .
The root nodes are represented by squares.
}
\label{fig.pn}
\end{figure}

The 1, 1, 2,  6 and 19 graphs for $p(n)$ up to
5 nodes are gathered in Figure \ref{fig.pn}.

The planted C-trees play a pivotal role in the recursive
build-up of C-trees. Given a multiset of one or more planted C-trees,
we may build a bundle of these by melting all their roots into
a single new root. Each node in a cycle of nodes of a C-tree
is such a bundled multiset.
The planted trees become branches of the C-tree,
see \cite[Fig 1]{LerouxASMQ16}.

\begin{remark}
In a generic setup, the empty graph would contribute $p(0)=1$,
but in our work one or more of
the root nodes of the planted C-trees
define a node in a cycle of a C-tree. We will be constructing
the C-trees by summing up the cases with $c\ge 1$ nodes
in the cycles such of explicit lengths.
The cycles must not be interrupted and at least one planted C-tree
must be fastened
at each node of a cycle (even if the C-tree has just a single node).
Therefore we
define $p(0)=0$ here.
\end{remark}

\begin{defn}
A planted C-forest is a multiset of planted C-trees where the
roots of the planted C-trees have been merged into a single root node.
\end{defn}
\begin{defn}
\begin{equation}
F(x) = \sum_{n\ge 1} f(n)x^n
=x+x^2+3x^3+9x^4\cdots
\end{equation}
is the generating function for planted C-forests with $n$ nodes.
\end{defn}

\begin{figure}
\includegraphics[scale=0.3]{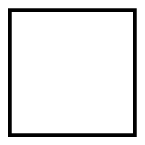}

\includegraphics[scale=0.3]{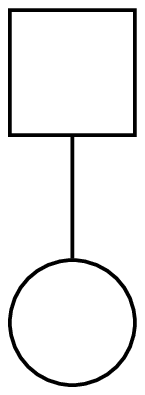}

\includegraphics[scale=0.3]{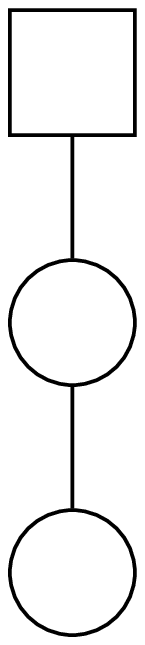}
\includegraphics[scale=0.3]{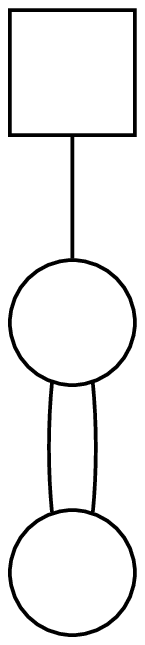}
\includegraphics[scale=0.3]{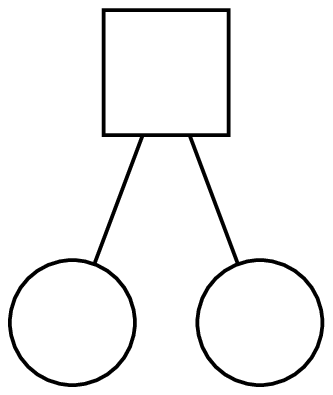}
\caption{The $f(1)=1$, $f(2)=1$ and $f(3)=3$ planted C-forests up to 3 nodes.}
\label{fig.fn}
\end{figure}

The first three terms of $F(x)$ are illustrated
in Figure \ref{fig.fn}.
Again, in this work the empty tree
is explicitly not represented, $f(0)=0$.

The construction of the planted C-forests means
removing the root of each planted tree, building a multiset of the trunks of 
these trees, and attaching them to a common root, so $F(x)/x$
is the Euler transform of $P(x)/x$ \cite{Flajolet},
explicitly (see e.g. \cite[eq 45]{HararyTAMS78})
\begin{equation}
F(x)/x = \prod_{n\ge 1}\frac{1}{(1-x^n)^{p(n+1)}} .
\label{eq.Feuler}
\end{equation}
This can also be written as \cite[(3.1.12)]{Harary}
\begin{equation}
F(x)/x = \exp\left[ \sum_{k\ge 1}\frac{1}{k}\left(
\frac{P(x^k)}{x^k}-1\right)
\right]
,
\end{equation}
where the generating function for the graphs with a marked entry node
(decapitated planted C-trees) is
\begin{equation}
C'(x) = \frac{P(x)}{x}-1 = \sum_{n\ge 1}p(n+1)x^n = x+2x^2+6x^3+19x^4+\cdots .
\end{equation}
\begin{remark}
The prime does not indicate differentiation but aligns the notation with Labelle, Robinson
and others.
\end{remark}

\begin{figure}
\includegraphics[scale=0.7]{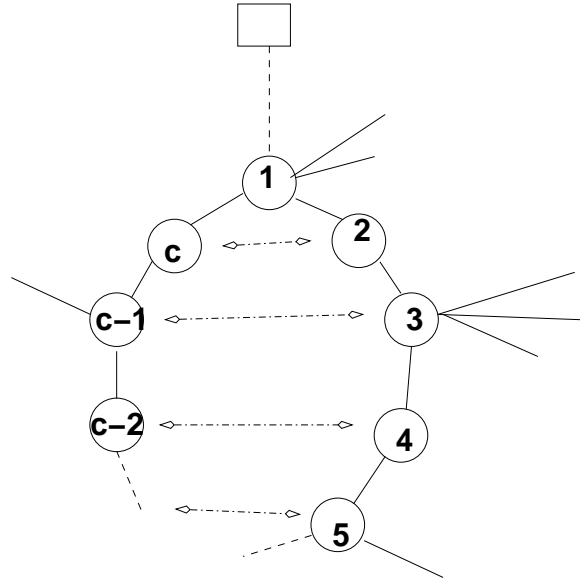}
\caption{Symmetry of the $S_2$ group implied in the 
cycles of planted C-trees. Nodes enumerated $1,2,\ldots,c$.}
\label{fig.wbrac}
\end{figure}

Computing $P(x)$  happens by summing over all possible
cycle length $c\ge 1$ that may occur in the cycle that is the neighbor
to the root node, and considering all multisets
of $c$ planted C-forests at the cycle nodes. The generic structure of
the cycles in planted C-trees is in Figure \ref{fig.wbrac}.

There is one unique node in the cycle (labeled 1 here, called the
entry node, which is also a cut node) that has smallest distance to the
root of the planted C-tree.
The planted C-forests attached to the nodes (possibly including
the entry node) are indicated by half-edges. Flipping the cycle
over while the entry node stays fixed, indicated
by the dashed lines with arrows, does not generate new graphs.
For fixed $c$,
the permutation group of order 2 of this flip has the following ``palindromic''
symmetry of swapping nodes:

\begin{tabular}{r|l}
$c$ & group generator (cycle representation) \\
\hline
1 & (1) \\
2 & (1)(2) \\
3 & (1)(23) \\
4 & (1)(3)(24) \\
5 & (1)(25)(34) \\
6 & (1)(4)(26)(35) \\
\end{tabular}

So the root of the cycle stays fixed, for even $c$
the opposite node stays also fixed, and nodes with
labels $i$ and $c+2-i$ swap places.
Polya's cycle index for (the automorphism group of) the
straight tree with $n$ nodes (the equivalent palindromic symmetry)
is known \cite[Table 2]{RobinsonJCT9}\cite[App.\ B]{MatharVixra1805}.
In our application the cycle has another fixed element where
the direction to the root of the planted tree enters the cycle, so
\begin{equation}
Z(S_2,c,t_1,t_2) = \left\{
\begin{array}{ll} 
1& c=0; \\
t_1(t_1^{c-1}+t_2^{(c-1)/2})/2,& $c=1,3,5,7,9,\ldots$; \\
t_1^2(t_1^{c-2}+t_2^{(c-2)/2})/2,& $c=2,4,6,8,\ldots$
\end{array}
 \right .
.
\end{equation}
The counting series at each node is $F(x)$, and Polya's method
of applying the $S_2$ symmetry to the generating function of the
planted C-trees yields
\begin{equation}
P(x) = x \sum_{c\ge 0} Z(S_2,c,t_1,t_2),\quad  t_i\to F(x^i).
\label{eq.polyaP}
\end{equation}
Summation over the geometric series in \eqref{eq.polyaP} yields
\begin{multline}
P(x)
= x\left(1
+\sum_{c=1,3,5,\ldots} Z(S_2,c,t_1,t_2)
+\sum_{c=2,4,6,\ldots} Z(S_2,c,t_1,t_2)
\right)
\\
= x\left(1
+\sum_{c=1,3,5,\ldots} t_1\frac{t_1^{c-1}+t_2^{(c-1)/2}}{2}
+\sum_{c=2,4,6,\ldots} t_1^2\frac{t_1^{c-2}+t_2^{(c-2)/2}}{2}
\right)
\\
=x\left(1+\frac{t_1(1+t_1)}{2} \left[\frac{1}{1-t_1^2}+\frac{1}{1-t_2}\right]\right)
=x\left(\frac12 +\frac{F(x)[1+F(x)]}{2[1-F(x^2)]}+\frac{1}{2[1-F(x)]}\right)
.
\end{multline}
On the right hand side we have substituted $1$ for the case
of an empty branch of the root node, ($c=0$), and multiplied by $x$
to complement the cycles by the root node.

Bootstrapping $P(x)$ starts from the rough estimator $P(x)\approx x$
including just the $n=1$ node,
derives $f(n)$ via \eqref{eq.Feuler} up to the
same $n$, inserts this into \eqref{eq.polyaP}
to compute $P(x)$ with an augmented truncation order $n+1$,
and cycles through this process:
\begin{multline}
p(n) = 
1, 1, 2, 6, 19, 67, 244, 934, 3665, 14755, 60466, 251690, 1060662, 4517568, 19413415,\\
    84073051, 366539371, 1607472753, 7086453177, 31385697280, 139586611475
\ldots,\quad n\ge 1;
\end{multline}
\begin{multline}
f(n) = 
1, 1, 3, 9, 31, 110, 417, 1617, 6466, 26335, 109109, 457968, 1944180, 8331081, \\
    35991543, 156581739, 685415080, 3016616752, 13340799273, 59254050302,
\ldots,\quad n\ge 1.
\label{eq.fnlist}
\end{multline}
If we exclude the graphs that have cycles of length 2, we get instead
\begin{multline}
P^{c\neq 2}(x)
= x\left(1
+\sum_{c=1,3,5,\ldots} Z(S_2,c,t_1,t_2)
+\sum_{c=4,6,\ldots} Z(S_2,c,t_1,t_2)
\right)
\\
=x\left(\frac12 -F^{c\neq 2}(x)^2+\frac{F^{c\neq 2}(x)[1+F^{c\neq 2}(x)]}{2[1-F^{c\neq 2}(x^2)]}+\frac{1}{2[1-F^{c\neq 2}(x)]}\right)
.
\end{multline}
The coefficients of the generating functions are
\begin{multline}
p^{c\neq 2}(n) = 
1, 1, 1, 3, 8, 24, 72, 231, 751, 2520, 8584, 29743, 104265, 369571,\\
    1321408, 4761876
\ldots,\quad n\ge 1;
\end{multline}
\begin{multline}
f^{c\neq 2}(n) = 
1, 1, 2, 5, 14, 41, 128, 410, 1356, 4576, 15723, 54767, 193062, 687203,\\
    2466837
\ldots,\quad n\ge 1
\end{multline}
which is illustrated by deleting the graphs with double edges
in Figures \ref{fig.pn} and \ref{fig.fn}.

\section{Rooted Skeleton Trees}
The planted C-trees and the planted C-forests employed above
are marking a single node. The root of a planted C-tree has
(at most) degree 1, whereas the degree of the root of a planted C-forest is not
bounded. A further type of rooted C-trees emerges
if the node of a skeleton tree is marked as a root
and that mark heaved to all nodes of the associated cycles
of the C-trees.
\begin{defn}
A skeleton-rooted C-tree is a C-tree where all nodes
of one of its cycles are marked.
\end{defn}
\begin{defn}
\begin{equation}
C^\cdot(x)=\sum_{n\ge 0}c^\cdot(n) x^n
= x+2x^2+5x^3+15x^4+49x^5+176x^6+239x^7+\cdots
\end{equation}
is the ordinary generating function for skeleton-rooted C-trees
with $n$ nodes.
\end{defn}

\begin{exa}
All ways of marking nonequivalent cycles of the C-trees with 
3 nodes of Figure \ref{fig.R3} lead to the $c^\cdot(3)=5$
graphs of Figure \ref{fig.Rdot3}.
\begin{figure}
\includegraphics[scale=0.3]{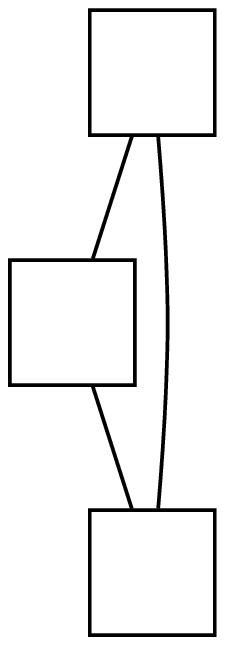}

\includegraphics[scale=0.3]{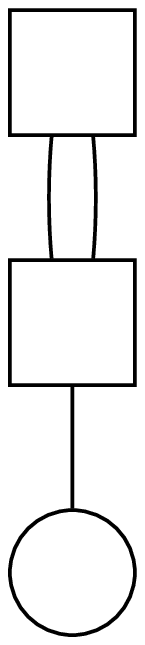}
\includegraphics[scale=0.3]{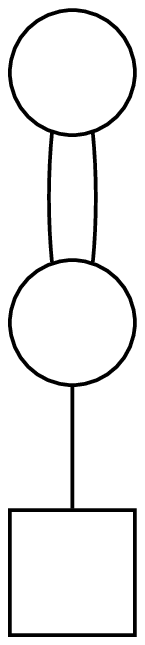}

\includegraphics[scale=0.3]{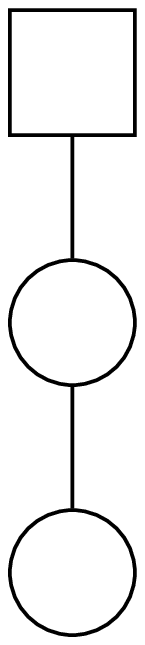}
\includegraphics[scale=0.3]{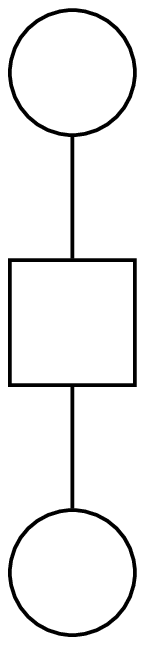}
\caption{The $c^{\cdot}(3)=5$ skeleton-rooted C-trees on 3 nodes.}
\label{fig.Rdot3}
\end{figure}
\end{exa}
\begin{exa}
All ways of marking noequivalent cycles of the C-trees with 
4 nodes of Figure \ref{fig.R4} lead to the $c^\cdot(4)=15$
graphs of Figure \ref{fig.Rdot4}.
\begin{figure}
\includegraphics[scale=0.3]{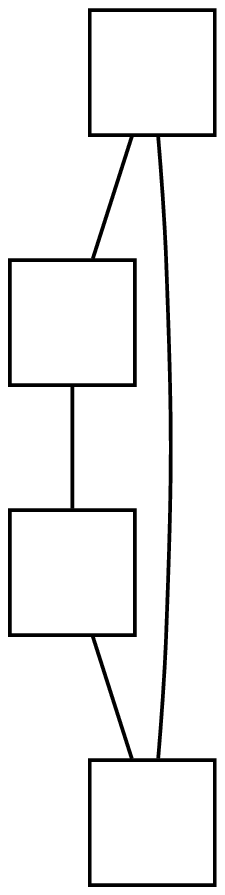}

\includegraphics[scale=0.3]{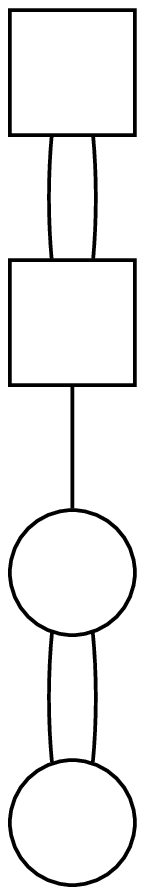}
\includegraphics[scale=0.3]{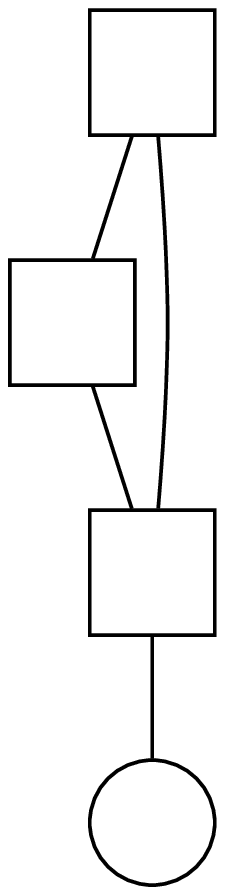}
\includegraphics[scale=0.3]{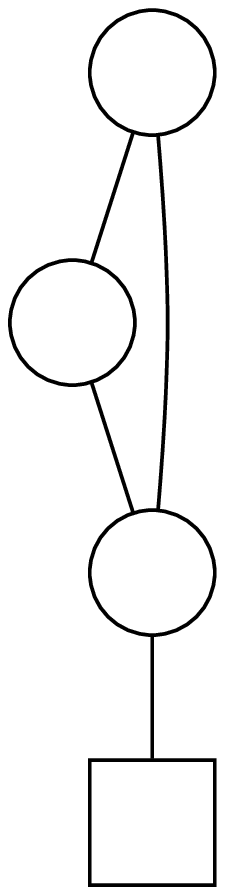}

\includegraphics[scale=0.3]{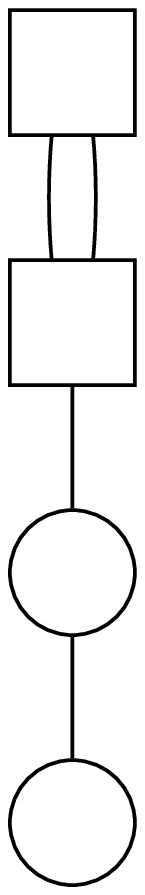}
\includegraphics[scale=0.3]{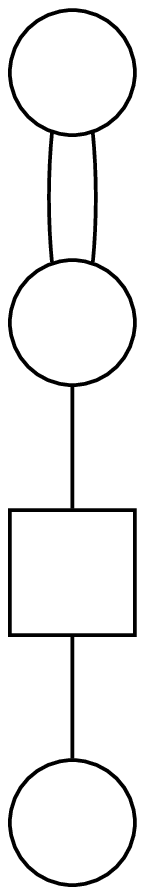}
\includegraphics[scale=0.3]{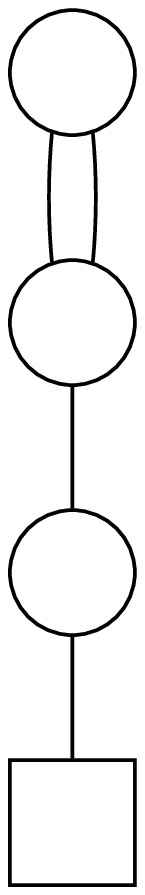}
\includegraphics[scale=0.3]{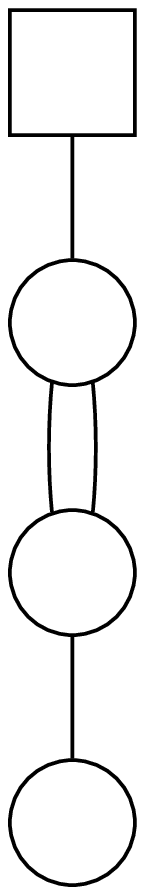}
\includegraphics[scale=0.3]{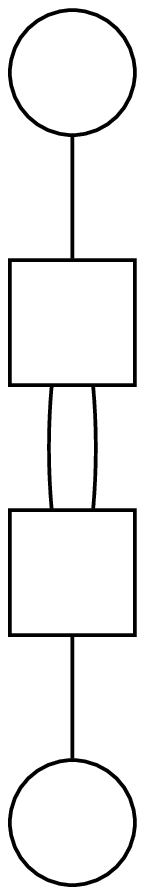}
\includegraphics[scale=0.3]{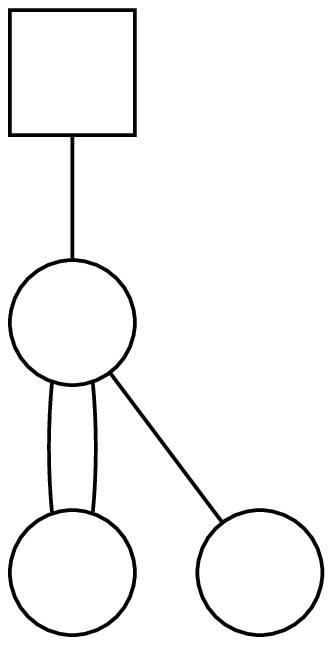}
\includegraphics[scale=0.3]{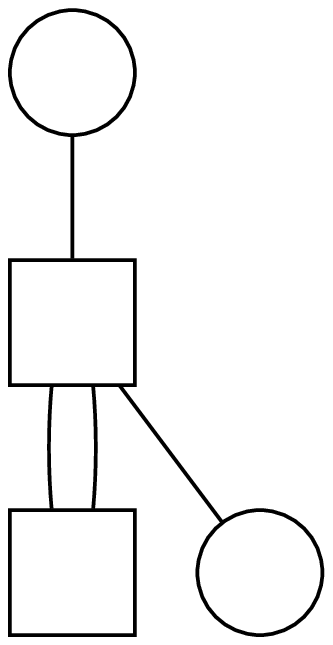}

\includegraphics[scale=0.3]{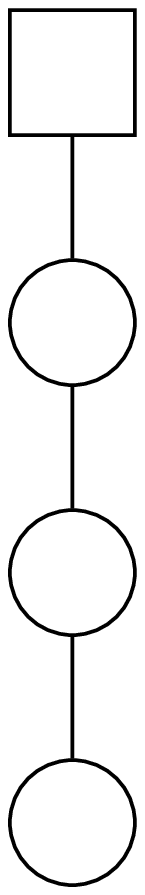}
\includegraphics[scale=0.3]{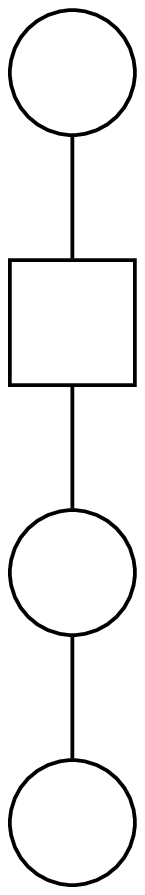}
\includegraphics[scale=0.3]{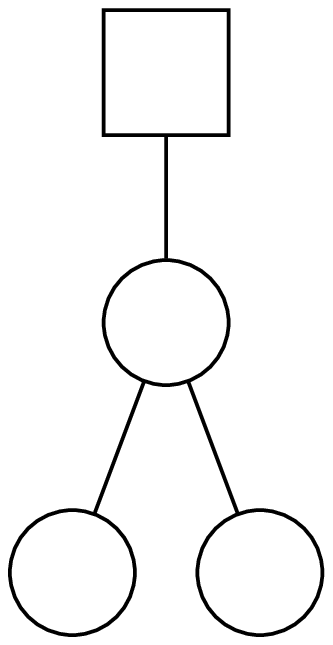}
\includegraphics[scale=0.3]{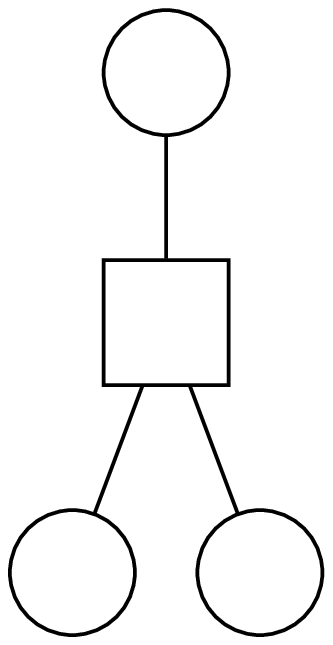}
\caption{The $c^{\cdot}(4)=15$ skeleton-rooted C-trees on 4 nodes derived from Figure \ref{fig.R4}.}
\label{fig.Rdot4}
\end{figure}
\end{exa}

\begin{figure}
\includegraphics[scale=0.7]{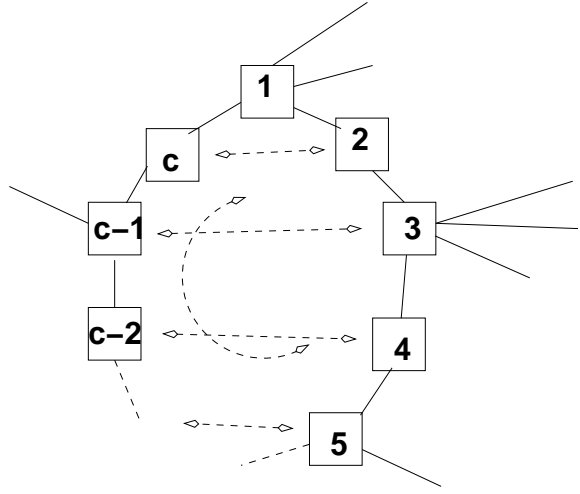}
\caption{Symmetry of the $D_c$ group implied in the 
root cycle of skeleton-rooted C-trees.}
\label{fig.rootC}
\end{figure}

The rooted cycle of a skeleton-rooted C-tree looks like Figure \ref{fig.rootC},
which is similar to Figure \ref{fig.wbrac} with the decisive difference
that there is no marked entry node; so the group of the symmetries
is no longer $S_2$ but the symmetry group of bracelets, the dihedral
group of $c$ elements.
The applicable cycle indices are the cycle indices of the Cyclic Groups \cite[(2.2.10)]{Harary},
\begin{equation}
Z(C_n,t_1,t_2,\ldots) = \frac{1}{n}\sum_{d\mid n} \varphi(d)t_d^{n/d}
\end{equation}
and the cycle indices of the Dihdral Groups \cite[(2.2.11)]{Harary}
\begin{equation}
Z(D_n,t_1,t_2,\ldots) =
\left\{
\begin{array}{ll}
\frac12 Z(C_n,t,t_1,t_2,\ldots) + \frac12 t_1t_2^{\lfloor n/2\rfloor},& n\, \textrm{odd} \\
\frac12 Z(C_n,t,t_1,t_2,\ldots) + \frac14 (t_1^2t_2^{n/2-1}+t_2^{n/2}) , & n\, \textrm{even}\\
\end{array}
\right.
.
\end{equation}
Each of the nodes in the rooted cycle is a planted C-forest, so
the generating function is
\begin{equation}
C^{\cdot}(x) = \sum_{c\ge 1} Z(D_c,t_1,t_2,\ldots),\quad t_i \to F(x^i).
\end{equation}
Insertion of \eqref{eq.fnlist} yields
\begin{multline}
c^\cdot(n)=
1, 2, 5, 15, 49, 176, 657, 2558, 10207, 41622, 172456, 724394, 3076455, \\
13189759,
    57004696, 248096112, 1086373375, 4782776966, \\
21157400729, 
93995763458, 419210682932,
    693975376003,
,\ldots \quad n\ge 1
.
\end{multline}
The variant of deleting graphs with cycles of length 2 is
\begin{equation}
C^{\cdot,c\neq 2}(x) = \sum_{c\ge 1,c\neq 2}
Z(D_c,t_1,t_2,\ldots),\quad t_i \to F^{c\neq 2}(x^i),
\end{equation}
\begin{multline}
c^{\cdot,c\neq 2}(n)=
1, 1, 3, 7, 19, 55, 168, 536, 1764, 5940, 20372, 70918, 249821, 888982, \\
3190384,
    11533780, 41962149, 153521353, 564448947, 2084469347, \\
7728325965, 6513869557
,\ldots \quad n\ge 1
.
\end{multline}

\section{Otter's Formula, Synthesis}

The counting series of the C-trees is synthesized from
the series of the planted C-trees and skeleton-rooted C-trees
with Otter's mapping between rooted trees and trees
\cite{OtterAM49,ClarkeQJMO10,HararyMMJ3,HararyAM101}\cite[(3.2.4)]{Harary} generalized
to enriched trees \cite[(0.13]{LerouxASMQ16}\cite[(2.4)]{LabelleJSC14}.

The generating function for pairs of these graphs (C-trees with
a marked edge that is a bridge) is
\begin{equation}
E_2(C'(x)) = \frac12[C'(x)^2+C'(x^2)] = x^2+2x^3+9x^4+31x^5+126x^6+492x^7+2014x^8+\cdots
\end{equation}
The generating function for unordered pairs of these graphs 
(C-trees with an oriented marked edge that is a bridge) is
\begin{equation}
[C'(x)]^2 = x^2+4x^3+16x^4+62x^5+246^6+984x^7+4009x^8+\cdots.
\end{equation}
Eventually the generating function $C(x)$ for the C-trees is derived via
\begin{equation}
C(x) + [C'(x)]^2 = C^{\cdot}(x)+E_2(C'(x)).
\end{equation}
The coefficients of the generating function are \cite[A317722]{EIS}
\begin{multline}
c(n) =
1, 2, 3, 8, 18, 56, 165, 563, 1937, 7086, 26396, 101383, 395821, 1573317,
6335511,\\
25825861,
    106344587, 441919711, 1851114466, 7809848543,
\ldots n\ge 1.
\end{multline}
The variant of not admitting cycles of length 2 in C-trees ---
which skips one C-tree of Figure \ref{fig.R2},
one C-tree of Figure \ref{fig.R3},
4 C-trees of Figure \ref{fig.R4},
10 C-trees of Figure \ref{fig.R5},
36 C-trees of Figures \ref{fig.R6a} and \ref{fig.R6b} ---
is
\begin{multline}
c^{c\neq 2}(n) =
1, 1, 2, 4, 8, 20, 48, 133, 374, 1124, 3439, 10923, 35245, 116128, 387729, \\
1312038,
4485906,
    15486546, 53900520, 188998450,\ldots
\quad
n\ge 1.
\end{multline}

\bibliographystyle{amsplain}
\bibliography{all}

\providecommand{\bysame}{\leavevmode\hbox to3em{\hrulefill}\thinspace}
\providecommand{\MR}{\relax\ifhmode\unskip\space\fi MR }
\providecommand{\MRhref}[2]{%
  \href{http://www.ams.org/mathscinet-getitem?mr=#1}{#2}
}
\providecommand{\href}[2]{#2}
\begin{thebibliography}{10}

\bibitem{ClarkeQJMO10}
L.~E. Clarke, \emph{On {Otter}'s formula for enumerating trees}, Q. J. Math.
  \textbf{2} (1959), no.~10, 43--45. \MR{0104595}

\bibitem{Flajolet}
Philippe Flajolet and Robert Sedgewick, \emph{Analytic combinatorics},
  Cambridge University Press, 2009. \MR{2483235}

\bibitem{HararyMMJ3}
Frank Harary, \emph{Note on the {P\'olya} and {Otter} formulas for enumerating
  trees}, Michigan Math. J. \textbf{3} (1955), 109--112. \MR{0078687}

\bibitem{HararyTAMS78}
\bysame, \emph{The number of linear, directed, rooted and connected graphs},
  Trans. Am. Math. Soc. \textbf{78} (1955), no.~2, 445--463. \MR{0068198}

\bibitem{HararyAnM58}
Frank Harary and Robert~Z. Norman, \emph{The dissimilarity characteristic of
  {Husimi} trees}, Ann. Math. \textbf{58} (1953), no.~1, 134--141. \MR{0055693}

\bibitem{Harary}
Frank Harary and Edgar~M. Palmer, \emph{Graphical enumeration}, Academic Press,
  New York, London, 1973. \MR{0357214}

\bibitem{HararyAM101}
Frank Harary and Geert Prins, \emph{The number of homeomorphically irreducible
  trees and other species}, Acta Math. \textbf{101} (1959), no.~1--2, 141--162.
  \MR{0101846}

\bibitem{HararyPNAS39}
Frank Harary and George~E. Uhlenbeck, \emph{On the number of {Husimi} trees},
  Proc. Natl. Acad. Sci. USA \textbf{39} (1953), no.~4, 315--322. \MR{0053893}

\bibitem{LabelleJSC14}
Gilbert Labelle, \emph{Counting asymmetric enriched trees}, J. Symbolic Comput.
  \textbf{14} (1992), no.~2--3, 211--242. \MR{1187233}

\bibitem{LerouxASMQ16}
Pierre Leroux and Brahim Miloudi, \emph{G\'en\'eralisations de la formule
  d'{Otter}}, Ann. Math. Qu\'ebec \textbf{16} (1992), no.~1, 53--80.
  \MR{1173933}

\bibitem{MatharVixra1805}
Richard~J. Mathar, \emph{Labeled trees with fixed node label sum},
  vixra:1805.0205 (2018).

\bibitem{OtterAM49}
Richard Otter, \emph{The number of trees}, Ann. Math. \textbf{49} (1948),
  no.~3, 583--599. \MR{0025715}

\bibitem{RobinsonJCT9}
Robert~W. Robinson, \emph{Enumeration of non-separable graphs}, J. Combin.
  Theory \textbf{9} (1970), no.~4, 327--356. \MR{0284380}

\bibitem{EIS}
Neil J.~A. Sloane, \emph{The {O}n-{L}ine {E}ncyclopedia {O}f {I}nteger
  {S}equences}, Notices Am.\ Math.\ Soc. \textbf{50} (2003), no.~8, 912--915,
  http://oeis.org/. \MR{1992789}

\end{thebibliography}

\end{document}